\RequirePackage{ifpdf}
\ifpdf 
\documentclass[pdftex]{sigma}
\else
\documentclass{sigma}
\fi

\begin{document}

\allowdisplaybreaks

\renewcommand{\thefootnote}{$\star$}

\renewcommand{\PaperNumber}{048}

\FirstPageHeading

\ShortArticleName{Singularities of Af\/f\/ine Schubert Varieties}

\ArticleName{Singularities of Af\/f\/ine Schubert Varieties\footnote{This paper is a
contribution to the Special Issue on Kac--Moody Algebras and Applications. The
full collection is available at
\href{http://www.emis.de/journals/SIGMA/Kac-Moody_algebras.html}{http://www.emis.de/journals/SIGMA/Kac-Moody{\_}algebras.html}}}

\Author{Jochen KUTTLER~$^\dag$ and Venkatramani LAKSHMIBAI~$^\ddag$}
\AuthorNameForHeading{J.~Kuttler and V.~Lakshmibai}

\Address{$^\dag$~Department of Mathematical {\rm \&} Statistical Sciences, University of Alberta, Edmonton, Canada}
\EmailD{\href{mailto:jochen.kuttler@ualberta.ca}{jochen.kuttler@ualberta.ca}}

\Address{$^\ddag$~Department of Mathematics, Northeastern University, Boston, USA}
\EmailD{\href{mailto:lakshmibai@neu.edu}{lakshmibai@neu.edu}}

\ArticleDates{Received September 11, 2008, in f\/inal form April 03,
2009; Published online April 18, 2009}

\Abstract{This paper studies the singularities of af\/f\/ine Schubert varieties
in the af\/f\/ine Grassmannian (of type $\mathrm{A}^{(1)}_\ell$). For
two classes of af\/f\/ine Schubert varieties, we determine the
singular loci; and for one class, we also determine explicitly the
tangent spaces at singular points. For a general af\/f\/ine Schubert
variety, we give partial results on the singular locus.}

\Classification{14M15; 14L35}

\newcounter{alphathm}

\newtheorem{thm}{Theorem}[section]
\newtheorem*{thmABC}{\stepcounter{alphathm}Theorem~\Alph{alphathm}}

\newtheorem{lem}[thm]{Lemma}
\newtheorem{prop}[thm]{Proposition}

\newtheorem{cor}[thm]{Corollary}
\theoremstyle{definition}
\newtheorem{defn}[thm]{Definition}
\newtheorem{rem}[thm]{Remark}

\newtheorem{ex}[thm]{Example}
\newtheorem{conj}[thm]{Conjecture}

\newcommand{\op}{\operatorname}

\newcommand{\fg}{{\mathfrak g}}
\newcommand{\bb}{{\mathfrak b}}
\newcommand{\PP}{{\mathbb P}}
\newcommand{\ZZ}{{\mathbb Z}}
\newcommand{\QQ}{{\mathbb Q}}

\newcommand{\Spec}{{\op{Spec\,}}}
\newcommand{\Mor}{{\op{Mor}}}
\newcommand{\Hom}{{\op{Hom}}}
\newcommand{\End}{{\op{End}}}
\newcommand{\Aut}{\op{Aut}}
\newcommand{\K}{K}
\newcommand{\CC}{{\mathbb C}}

\newcommand{\Lie}{\op{Lie}}

\newcommand{\pp}{{\mathfrak p}}
\newcommand{\uu}{{\mathfrak u}}
\newcommand{\uum}{{\mathfrak u^-}}
\newcommand{\phit}{{\Phi_\tau^+}}
\newcommand{\phitm}{{\Phi_\tau^-}}

\newcommand{\Proj}{\op{Proj}}
\newcommand{\bx}{\mathbf{x}}

\newcommand{\calG}{{\mathcal G}}
\newcommand{\GL}{{\op{GL}}}
\newcommand{\SL}{{\op{SL}}}
\newcommand{\calB}{{\mathcal B}}
\newcommand{\calP}{{\mathcal P}}
\newcommand{\calF}{{\mathcal F}}
\newcommand{\calW}{{\mathcal W}}
\newcommand{\ord}{{\op{ord}}}

\newcommand{\Gr}{{\op{Gr}}}
\newcommand{\IGr}{{\op{Gr}_\infty(\K)}}
\newcommand{\real}{{\op{re}}}
\newcommand{\im}{{\op{im}}}

\newcommand{\Ga}{{\mathbb G_a}}
\newcommand{\Gm}{{\mathbb G_m}}
\newcommand{\calX}{{\mathcal X}}
\newcommand{\reg}[1]{{\mathcal O({#1})}}
\newcommand{\Sn}{{\mathfrak S_n}}

\newcommand{\tor}{{\op{Tor}}}

\newcommand{\hh}{{\mathfrak h}}
\newcommand{\CM}{\Delta}
\newcommand{\YW}{{\op{YW}}}
\newcommand{\Sc}{{\mathfrak S}}

\newcommand{\calH}{{\mathcal H}}

\newcommand{\ws}{{w^s}}
\newcommand{\wsn}{{w^s_0}}
\newcommand{\sing}{{\varphi}}

\newcommand{\name}[1]{{\textsc{#1}}}

\newcommand{\codim}{{\op{codim}}}

\newcommand{\supp}{{\op{supp}}}

\newcommand{\BB}{{\mathfrak B}}

\newcommand{\hPhi}{{\widehat \Phi}}
\newcommand{\halpha}{{\hat \alpha}}
\newcommand{\hbeta}{{\hat \beta}}
\newcommand{\hT}{{\widehat T}}
\newcommand{\hW}{{\widehat W}}

\newcommand{\TV}{{\mathfrak T}}
\newcommand{\BV}{{\mathfrak B}}
\newcommand{\WV}{{\mathfrak W}}
\newcommand{\PhiV}{{\Phi(V)}}

\newcommand{\ev}{\op{ev}}

\newcommand{\tu}[1]{\textup{#1}}

\newcommand{\da}{\downarrow}
\newcommand{\up}{\uparrow}

\newcommand{\vdim}{{\op{vdim}}}

\section{Introduction}

Schubert varieties are important objects in the theory of
algebraic groups, representation theory, and combinatorics. The
determination of their singularities is a classical problem and
has been studied by many authors. Very conclusive results are
available for groups of type $A$ (see e.g.\ \cite{LSe, LS,BW, CO, KLR,MA}; and for arbitrary
types \cite{BP,KUMI,CK}). In this note we
investigate singularity properties of the natural generalization
of Schubert varieties to \emph{affine Schubert varieties}. In the
af\/f\/ine setting, the question has not been settled yet, although,
there are certainly results available (in particular, \cite{KUMI}
applies as well). Recently, af\/f\/ine Schubert varieties (in all
types) have been studied by several authors
(cf.~\cite{BM,em,ju,mov}). While in~\cite{em,ju,mov}, the authors
study the singularities of $\mathcal{P}$-stable af\/f\/ine Schubert
varieties (see Section~\ref{subsec:PSTABLE}),  in~\cite{BM}, the authors classify the
smooth and rationally smooth Schubert varieties.

The most classical Schubert varieties are  the Schubert varieties
in Grassmannians, and the f\/irst generalization is therefore to the
af\/f\/ine Schubert varieties in the af\/f\/ine Grassmannian of type~$A^{(1)}$. So let us f\/ix an algebraically closed f\/ield $\K$ of
characteristic zero, and denote by \mbox{$A = \K[[t]]$} the ring of
formal power series with quotient f\/ield $F = \K((t))$, the ring of formal Laurent series. Then $\SL_n(A)$
and $\SL_n(F)$ both are the $\K$-points of ind-varieties over~$\K$, denoted by $\calP$ and $\calG$, and $\calP \subset \calG$.
The \emph{affine Grassmannian} is then the quotient ind-variety~$\calG / \calP$. Mi\-micking the classical situation the af\/f\/ine
Schubert varieties are the $\calB$-orbit closures in~$\calG/\calP$, where $\calB \subset \calP$ is the subgroup of
elements where the (strictly) upper triangular entries are
divisible by~$t$; more formally $\calB = \ev^{-1}(B)$, where $B
\subset \SL_n(\K)$ is the Borel subgroup of lower triangular
matrices and $\ev \colon \calP \to \SL_n(\K)$ is the evaluation
homomorphism sending $[g_{ij}(t)]$ to~$[g_{ij}(0)]$. Let~$T$ be
the maximal torus consisting of diagonal matrices in $\SL_n(\K)
\subset \calP$, and let $S = K^*$ be the one-dimensional torus in
$\Aut(\calG)$ coming from the action of $S$ on $F$ by rotating the
loops, i.e.\ $s \in S$ sends $g(t)$ to $g(st)$. As the $S$ and
$T$-actions commute, putting $\hT = T \times S$ we obtain an
$n$-dimensional torus acting on $\calG$, $\calP$ and $\calB$, and
therefore on any af\/f\/ine Schubert variety. Each $\calB$-orbit
contains a unique $\hT$-f\/ixed point; the $\hT$-f\/ixed points in
$\calG / \calP$ are parameterized by~$\hW^\calP$, a set of
representatives of~$\hW/W$,
 $\hW$ (resp.~$W$) being the af\/f\/ine Weyl group (resp. the Weyl group)
 of type $\mathrm{A}^{(1)}_{n-1}$ (resp.\ $\mathrm{A}_{n-1}$). In fact,
$\hW^\calP$ has a natural identif\/ication with $\ZZ^{n-1}$ embedded
in $\ZZ^n$ as the sublattice consisting of points in $\ZZ^n$  with
coordinate sum being equal to zero. Let $\ge$ denote the partial
order on $\hW^\calP$ induced by the partial order on (the Coxeter
group) $\hW$ (with respect to the set of simple roots associated
to $\mathcal{B}$). For $w \in \hW^\calP$, let
\[ X(w):={\underset{\{v\in\hW^\calP,\,v\le w\}}{\bigcup}}\calB v
\]
be the \emph{affine Schubert variety} in $\calG / \calP$  associated to
$w$. (Thus for $v,w\in\hW^\calP$, we have, $v \leq w$ if and only
if $X(v) \subseteq X(w)$).

In this paper, for studying the af\/f\/ine Grassmannian and the af\/f\/ine
Schubert varieties, we make use of  a canonical embedding of
af\/f\/ine Grassmannian into the \emph{infinite Grassmannian}~$\Gr(\infty)$ over~$\K$ (cf.\ Section~\ref{sec:PRELIMINARIES}). We
brief\/ly explain below our approach.

 $\Gr(\infty)$ being the
inductive limit of all f\/inite dimensional Grassmannians, we obtain
a canoni\-cal identif\/ication of an af\/f\/ine Schubert variety $X(w)$ as
a closed subvariety of a suitable Grassmannian $G(d,V)$ (the set
of $d$-planes in the vector space $V$), in fact, as a closed
subvariety of a suitable classical Schubert variety in $G(d,V)$.
Further, as a subset of $\Gr(\infty)$, we get an identif\/ication of
$\calG / \calP$ with the set of $A$-lattices in $F^n$ (i.e., free
$A$-submodules of $F^n$ of rank~$n$). For instance, the element
$(c_1,\dots,c_n)\in\hW^\calP$ corresponds to the $A$-span of
$\{t^{c_1}e_1,\dots,t^{c_n}e_n\}$ (here,$\{e_1,\dots,e_n\}$ is
the standard $F$-basis for $F^n$). Given
$w=(c_1,\dots,c_n)\in\hW^\calP$, there exists an $s > 0$ such
that $w \leq \wsn$ with $\wsn :=
(-s(n-1),s,s,\dots,s)(\in\hW^\calP)$, and hence any $X(w)$ may be
thought of as a subvariety of $X(\wsn)$, for a suitable $s$. We
have an identif\/ication (cf.\ Section~\ref{iden})
\[
X(\wsn)=\{A{\texttt{-}}{\mathrm{lattice\, }}L|\;\; t^{-s(n-1)}L_0
\supset L\supset t^sL_0,\,\mathrm{and} \,\dim_K(L/t^sL_0)=sn\},
\]
$L_0$ being the standard lattice (namely, the $A$-span of
$\{e_1,\dots,e_n\}$). Setting $V\!:=t^{-s(n-1)}L_0/t^{s}L_0$, we
have $\dim V=sn^2$; further, the map $f_s: X(\wsn)\rightarrow
G(d,V)$, $L\mapsto L/t^sL_0$ (where \mbox{$d:=sn$}) identif\/ies $X(\wsn)$ as
a closed subvariety of $G(d,V)$. Denoting $u:=1+t$, the unipotent
endomorphism of $V$, $v\mapsto v+tv$, we have that $u$ induces an
automorphism of $G(d,V)$ and $f_s$ identif\/ies $X(\wsn)$ with
$G(d,V)^u$ (the f\/ixed point set of $u$ with the reduced scheme
structure) (cf.\ Proposition~\ref{kap0}).  Moreover,  for each
af\/f\/ine Schubert variety $X(w)\subseteq X(\wsn)$, we have that $w
\in G(d,V)$ is a $\TV$-f\/ixed point (for a suitable maximal torus
$\TV$ of $\GL(V)$), giving rise to a classical Schubert variety
$Y(w)$ (with respect to a suitable Borel subgroup $\BV$) which is
$u$-stable, and we have an identif\/ication: $X(w)= Y(w)^u$. Thus we
exploit this situation to
 deduce properties for af\/f\/ine Schubert varieties. In particular,
 for two classes of af\/f\/ine Schubert varieties contained in
$X(\wsn)$, we determine explicitly the singular loci; further, for
one of the two classes, we also determine the tangent spaces at
singular points. In order to describe our results, given $w =
(c_1,c_2,\dots,c_n) (\in \hW^\calP \subset \ZZ^n)$, such that
$X(w)\subseteq X(\wsn) (\subset G(d,V))$, we def\/ine  $L(w): =
(l_1,l_2,\dots,l_{2n})$, where
\[
l_i =
\begin{cases}s-c_i,&
\mathrm{if\ } i \leq n,\\
  s-c_{i-n} + 1,&\mathrm{if\ } i > n.
  \end{cases}
  \]

The f\/irst class of af\/f\/ine Schubert varieties that we consider consists of
the $\calP$-stable af\/f\/ine Schubert varieties.
Note that the $\calP$-stable af\/f\/ine Schubert varieties $X(w)$'s contained in
$X(\wsn)$ may be characterized by the corresponding $L(w)$'s:
$X(w)$ is $\calP$-stable if and only if $l_1\ge l_2\ge\cdots\ge
l_n$. The
Schubert variety $X(\wsn)$ is an example of a $\calP$-stable
Schubert variety, and an important at that, in view of its
relationship (cf.~\cite{LU2}) with nilpotent orbit closures in
$\Lie(G)$ (for the adjoint action of $G$ on $\Lie(G)$ $G$ being
$\GL_n(K)$). To be very precise, for $s=1$, $X(\wsn)$ contains the
variety of nilpotent matrices as an open subset; moreover, we have
a bijection between \{nilpotent orbit closures\} and
\{$\calP$-stable af\/f\/ine Schubert subvarieties of $X(w_0^1)$\} -- by
Lusztig's isomorphism (cf.~\cite{LU2}), a nilpotent orbit closure
gets identif\/ied with an open subset (namely, the ``opposite cell'')
of an unique Schubert subvariety of $X(w_0^1)$. For a
$\calP$-stable af\/f\/ine Schubert variety $X$, we determine
Sing$\,X$, the singular locus of $X$; of course, this result was
f\/irst proved by Evens and Mircovi\'c (cf.~\cite{em}). We give yet
another proof of their result, which is that the regular locus of a $\calP$-stable Schubert variety is precisely the open $\calP$-orbit. We also determine explicitly (cf.\
Corollary~\ref{cor:CURVE_COR}, Theorems~\ref{kap1},~\ref{thm:ONE_STRING_TANGENT}) $T_xX(\wsn)$, the tangent
space to $X(\wsn)$ at $x\le \wsn$; our description is in terms of
$T_xY(\wsn)$, for any $s>0$.

The second class of af\/f\/ine Schubert varieties consists of $X(w)$'s
such that $L(w)$ admits two indices $i\leq n$, $i \leq j < i + n$
, such that for $k = i,\dots,i+n-1,k\ne j$, we have $l_k$ is
independent of $k$, and less than or equal $l_j$ (cf.\
Section~\ref{sec:ONESTRING}); we say that $w$ consists \emph{of
one string}. Note that $\wsn$ consists of one string. For these
$X(w)$'s, we show (cf. Theorem \ref{thm:ONE_STRING_TANGENT}):

\begin{thmABC}\label{theoremA}
Let $w$ consist of one string. Then
$T_xX(w)=T_xX(\wsn)\cap T_xY(w)$.
\end{thmABC}

We also prove the rational smoothness of $X(\wsn)$ and certain of
the Schubert varieties of the ``one-string'' type (cf.\ Theorem~\ref{thm:RATIONALLY_SMOOTH} and its corollary), which is also obtained in~\cite{BM}.

The realization of the af\/f\/ine Schubert variety $X(w)$ as a closed
subvariety of the classical Schubert variety $Y(w)$ ($\subset
G(d,V)$) enables us to construct certain singularities as
explained below.
We say $P = (i,j)$ ($i < j$) is an \emph{imaginary pattern} in
$L(w)$ if $l_i > l_{j} + 1$. We may assume $i \leq n$.

\begin{thmABC}\label{theoremB}
Let $P = (i,j)$ be an imaginary pattern in $w$ and let $w_P$ be
defined by $L(w_P)_{i} = l_i - 1$, $L(w_P)_j = l_j  + 1$, and
$L(w_P)_k=L(w)_k$, $k\ne i,j$. Then $w_P$ is a singular point of
$X(w)$.
\end{thmABC}

It turns out that Theorem~B is enough to describe the maximal
singularities of the two classes of Schubert varieties described
above.

The reason why $w_P$ is singular is simply that $T_{w_P}(X(w))$
contains a tangent line whose $\hT$-weight is an imaginary root.
As these are never tangent to $\hT$-stable curves, $w_P$ has to be
singular. Of course, another possible reason for singularity is
that there may be too many of such curves. It is not hard to
construct points in most Schubert varieties where this is the
case: Let $P\colon i < g < j < k \leq 2n$ be a sequence with $i
\leq n$, such that $j < i + n$, $k < g + n$, and $l_i \geq l_j >
l_g \geq l_k$. We call $P$ a \emph{real pattern of the first
kind}. Def\/ine $w_P$ by putting $(l_g,l_k,l_i,l_j)$ in~$L(w)$ at
positions $(i,g,j,k)$ (see Section~\ref{sec:PATTERNS}). Similarly,
let $Q \colon i < j < g < k \leq 2n$ be a~sequence of integers
such that $i \leq n$, $g < i + n$, $k < j + n$, and $l_j
> l_i \geq l_k > l_g$. We refer to~$Q$ as a \emph{real
pattern of the second kind}. Def\/ine $w_Q$ by putting
$(l_g,l_k,l_i,l_j)$ in $L(w)$
 at positions $(i,j,g,k)$.

\begin{thmABC}
 If $w$ admits a real pattern $P$ of any kind then
$w_P$ is a singular point in $X(w)$.
\end{thmABC}

As in the classical setting, the geometric explanation why $w_P$ is singular is that the dimension of $T_x(X(w))$ is too big due to the presence of too many $\hT$-invariant curves each of whom contributes a line in $T_x(X(w))$.

\begin{rem}
We observe that the relative order of the lengths
$l_i$, $l_j$, $l_g$, $l_k$ in both of the real patterns is almost the same
as for the Type I and II patterns for classical Schubert varieties
(cf.~\cite{LS}). The dif\/ference is that we allow non-strict
inequalities at some places. Also, if $P$ is such a pattern the
relative order of the lengths in~$w_P$ is the same as in the
singularity constructed from the pattern in the classical setting.

Of course this begs the question whether the results of the classical setting could be applied directly to show that $w_P$ is singular. While in some examples this seems indeed possible we haven't been able so far to make this precise except for ``obvious'' cases. Hopefully we will be able to address this question more satisfactorily in some future work. We give some indication in Remark~\ref{rem:PATTERN}.

The same can be said for the more general question, whether we can actually formulate
these pattern in terms of the Weyl group elements and relate it to work of Billey--Braden~\cite{BB}
or Billey--Postnikov~\cite{BP} in the f\/inite case.
\end{rem}

Of course, this raises the question whether all (maximal)
singularities arise in this fashion. So far we haven't been able
to answer this question. However, in many examples it is true, if
one allows two degenerated cases of real patterns as well (as
discussed in Section~\ref{sec:PATTERNS}).

The paper is organized as follows: In
Section~\ref{sec:PRELIMINARIES}, we establish the basic notation
and conventions used throughout the text, and we will collect some
elementary results which will help further on.
Section~\ref{sec:REFLECTIONS} introduces the main combinatorial
tools as well as the notion of ``small'' ref\/lections and the
language of up-/down-exchanges, which serve as a tool in
describing the singularities later. In Section~\ref{sec:RELATION}
we investigate the relation between a Schubert variety $X(w)$ and
the classical Schubert variety $Y(w)$ which contains it, as far as
tangent spaces are concerned, and we introduce \emph{imaginary}
and \emph{real} tangents.  Section~\ref{sec:PATTERNS} again
returns to combinatorics, precisely def\/ining the various patterns
in $w$ which give rise to singularities, and proving Theorems~A
and B. The f\/inal Section~\ref{sec:ONESTRING} then applies these
results to two classes of Schubert varieties, those consisting of
``one string'', and those that are $\calP$-stable.

\section{Preliminaries}\label{sec:PRELIMINARIES} As explained in the
Introduction, for the study of the af\/f\/ine Schubert varieties in
$\calG/\calP$, we make use of a canonical embedding (as an
Ind-subvariety) of $\calG/\calP$ into $\Gr(\infty)$, the inf\/inite
Grassmannian. We shall now describe this embedding. For details,
we refer the readers to~\cite{KUM} and~\cite{MAG}.

\subsection[The affine and infinite Grassmannians]{The af\/f\/ine and inf\/inite Grassmannians}

We will keep the notation already established in the Introduction.
Consider the $\K$-vector space~$\K^\infty$ of $\K$-valued
functions on $\ZZ$ that vanish on ``very negative'' values, i.e.\
$\K^\infty = \{ f \colon \ZZ \to \K \mid f(i) = 0; \, i \ll 0\}$. For
each $i\in\mathbb{Z}$, there exists a canonical element $e_i\in
K^\infty$ def\/ined by $e_i(j) = \delta_{ij}$. Then every $f \in
\K^\infty$ may be written formally as $f = \sum_{i\in\ZZ} f_i e_i$
where $f_i = f(i)$.
Let $E_r = \{f \in \K^\infty \mid
f(i) = 0, \forall \,i < r\}$ be the ``span'' of $e_r,e_{r+1},\dots$.
The \emph{infinite Grassmannian $\Gr(\infty)$} over $\K$ is by
def\/inition the Ind-scheme obtained as the direct limit of usual Grassmannians as outlined below. Its $\K$-valued points are given by
linear subspaces $E$ of $\K^\infty$, such that for some $r
> 0$, $E_r \subset E$  and $E/E_r$ is
f\/inite-dimensional. Obviously, for such an $E$ we may increase~$r$
if necessary so that $E_{r+1} \subset E \subset E_{-r}$. As both
$E$ and $E_1$ contain $E_{r+1}$, we have
\[
\dim E/E_{r+1} = r + \dim E/(E_1 \cap E) - \dim E_1 /(E_1 \cap
E).
\] We set $\vdim(E) := \dim E_1 / (E_1 \cap E) - \dim E/(E_1
\cap E)$. For any $\K$-vector space $V$ and any positive integer
$d$, let $G(d,V)$ denote the Grassmannian of $d$-planes in $V$.
Then, for $s$ suf\/f\/iciently large, we have that $E$ is naturally an
element of $G(r-\vdim(E),E_{-s}/E_{r+1})$. The Grassmannians
$\{G(r-i$, $E_{-s}/E_{r+1}))_{s>0,r\geq i}\}$ ($i$ being f\/ixed) form
a direct system of varieties, with $G(r-i,E_{-s}/E_{r+1})$ and
$G(r'-i,E_{-s'}/E_{r'+1})$ both mapping naturally to
$G(r+r'-i,E_{-s-s'}/E_{r+r'+1})$. Its limit Ind-variety is denoted
$\Gr(\infty)_i$, and it parameterizes all $E$ with $\vdim (E) =
i$. Thus, $\Gr(\infty) = \bigcup_{i \in \ZZ} \Gr(\infty)_i$
carries a natural Ind-variety structure.

Let $\mathfrak X$ be the set of all $A$-lattices in $F^n$. Then
$\mathfrak X$ is naturally an algebraic subset of the set of
$K$-valued points of $\Gr(\infty)$ and therefore carries a
structure as an Ind-scheme. Indeed, let $F^n \to \K^\infty$ be the
isomorphism that sends $t^iv_j$ to $e_{j + in}$, $\{v_j,1\le j\le
n\}$ being a $F$-basis of $F^n$. Since for every $A$-lattice $L
\subset F^n$, $F^n /L$ is a torsion module, $L$ contains $t^rA^n$
for a suitably large $r
> 0$. On the other hand, if $\{u_i,1\le j\le n\}$ is an $A$-basis
for $L$, and $\ord\, u_i \geq N$, $\forall\, i$ for some $N \in \ZZ$,
then $L \subset t^NA^n$. Thus the image of $L$ in $\K^\infty$
def\/ines naturally an element of $\Gr(\infty)$; further, we have,
\[
\mathfrak X=\{E\in\Gr(\infty)|tE\subseteq E\}.
\] Here, via
 the identif\/ication of $F^n$ with
$\K^\infty$, $t$ acts on $\K^\infty$  as $te_i = e_{i+n}$; in the
sequel, we shall denote the map
$\K^\infty\rightarrow\K^\infty,\,e_i \mapsto e_{i+n}$, by $\tau$.
Notice that $\tau$ def\/ines a nilpotent map on each
$E_{-r}/E_{r+1}$, while the left-multiplication by $t$ def\/ines an
automorphism of $\Gr(\infty)$.

To describe the Ind-variety structure on $\calG/\calP$,  we
consider the natural transitive action of $\GL_n(F)$ on $\mathfrak
X$; observe that $\GL_n(A)$ is the stabilizer of the standard
lattice $A^n \subset F^n$. Hence, we get an identif\/ication
$\GL_n(F)/\GL_n(A)\cong\mathfrak X (\subset \Gr(\infty))$. Now
$\calG$ acts naturally on $\Gr(\infty)$ and on $\mathfrak X$, and
it is easily seen that $\mathfrak X \cap \Gr(\infty)_0$ is exactly
one $\calG$-orbit with $\calP$ as the stabilizer of the standard
lattice. Thus we obtain an embedding
$\calG/\calP\hookrightarrow\Gr(\infty)_0$ identifying
$\calG/\calP$ as an Ind-subvariety of $\Gr(\infty)_0$.

\subsection{The Weyl group}\label{weyl}
Let $G = \SL_n(K)$ with Lie algebra $\fg$.
Let $\Phi$ be the root system of $(G,T)$, and let $W = \Sn$ be
the associated Weyl group. The associated af\/f\/ine root system
$\hPhi$ then is by def\/inition the set of roots of $\hT$ in $\fg
\otimes_\K \K[t,t^{-1}]$. It may be identif\/ied with $\Phi \times
\ZZ \cup \{ 0 \} \times \ZZ$, and we write $\delta$ for $(0,1) \in \hPhi$. The elements of
$\ZZ\delta$ are called \emph{imaginary roots}, and all other
elements are called \emph{real roots}. For $\halpha = \alpha +
h\delta$ with $\alpha \in \Phi$ we put $\Re(\halpha) =\alpha$. A root $\halpha = \alpha + h\delta \in \hPhi$ is \emph{positive}  if $h > 0$, or $h = 0$ and $\alpha$ is positive (in the usual sense with respect to $\Phi$); otherwise, $\halpha$ is called \emph{negative}.

Let
$\Sc_\infty$ be the group of permutations of $\ZZ$, and $\tau$ the
element in $\Sc_\infty$:  $\tau(i) = i + n$, $i \in \ZZ$. Let
$\widetilde W = \{ \sigma \in \Sc_\infty \mid \tau \sigma = \sigma
\tau\}$. Clearly $\widetilde W = W \rtimes \ZZ^n$, where $W$
embeds naturally into $\Sc_\infty$, acting on the intervals
$[1+kn,(k+1)n]$, and $\ZZ^n$ embeds as $c = (c_1,\dots,c_n)$ maps
to $\tau^c$ with $\tau^c(i + kn) := i + (k+c_i)n$. The Weyl group
$\hW$ of $(\calG,\hT)$ may be naturally identif\/ied with a~subgroup
of $\widetilde W $.

$\hW \subset \widetilde W$ is given as the set of those $(w,c) \in
\widetilde W$ with $\sum_i c_i = 0$.  It is generated by
ref\/lections~$s_\halpha$ associated to the real roots $\halpha \in
\Phi$: For a root $\alpha = (ij) \in \Phi$ write $c_\alpha = e_j -
e_i \in \ZZ^n$ where $\{e_k\}$ denotes the standard basis of
$\ZZ^n$.  If $\halpha = \alpha + h \delta \in \hPhi$ with $\alpha$
being a positive root in $\Phi$, then $s_\halpha =
(s_\alpha,hc_\alpha) \in \widetilde W$, where $s_\alpha \in W$ is
the permutation associated to $\alpha$. If $\alpha = (ij)$ with $i
> j$, then
\[s_\halpha(q)=\begin{cases}j + (k-h)n,&\mathrm{\ if\ }q=i+kn,\\
i + (k+h)n,&\mathrm{\ if\ }q=j+kn,\\
q,&\mathrm{\ if\ }q\not\equiv i\mathrm{\ or\ }j \ (\mathrm{mod\,}n).
\end{cases}\]

Further, $s_{-\halpha}$ is then def\/ined as $s_\halpha$.

\subsection{Schubert varieties}
As mentioned in the introduction, each $\calB$-orbit on
$\calG/\calP$ contains a unique $\hT$-f\/ixed point. In fact, these
f\/ixed points form one orbit under the natural action of $\hW$.
These are best described in the language of lattices. Clearly an
$A$-lattice $L \subset F^n = \K^\infty$ is normalized by $\hT$ if
and only if it has a basis of the form
$e_{i_1},e_{i_2},\dots,e_{i_n}$; it is clear that as an element of
$\Gr(\infty)$, $L$ has the form $L = V_0 \oplus E_{r}$ for some
$r>0$, and some subspace $V_0$ spanned by a subset of the natural
basis of $E_{-r}/E_{r}$. Therefore $L$ is uniquely determined by
the ascending sequence $w(L)$ of integers, describing which $e_i$
lie in $L$: $w(L) = (w_1,w_2,\dots)$, and $w_i$ occurs if and only
if $e_{w_i} \in L$. Notice that eventually $w(L)$ agrees with the
natural sequence since $E_r \subset L$ for some $r$. One checks
easily that $L = \tau^cE_1$ for a suitable $c$, and in fact
\begin{equation}\label{eq:WT}
\ZZ^n \to \mathfrak X^\hT, \qquad c \mapsto \tau^cE_1
\end{equation}
is a bijection. Notice that $\vdim (\tau^c E_1) = \sum_i c_i$.
Further, as seen above, $\mathfrak X \cap \Gr(\infty)_0 =
\calG/\calP$, and the $\hT$-f\/ixed points therein are given by the
image of $\ZZ^{n-1}$ under the map (\ref{eq:WT}). We thus obtain a
notion of Schubert variety in all of $\mathfrak X$ (these are of
course the Schubert varieties for $\GL_n(F)/\GL_n(A)$). Recall
that $t$ acts as an automorphism of $\Gr(\infty)$ and $\mathfrak
X$, which commutes with the action of $\calG$. Clearly
$t(\Gr(\infty)_i) = \Gr(\infty)_{i+n}$ for all $i$, and under this
map $\calB$-orbit closures are sent to $\calB$-orbit closures. For
any $w \in \ZZ^n$ we put $X(w) = \overline{\calB w} \subset
\mathfrak X$. Here we are only interested in those $X(w)$ which
lie in $\Gr(\infty)_0$.

\subsection[The Schubert variety $X(\wsn)$]{The Schubert variety $\boldsymbol{X(\wsn)}$}\label{iden}

For $s > 0$, let $\wsn = \tau^{(-s(n-1),s,\dots,s)}$. Clearly, the
lattice $L_{\wsn}(=\tau^{(-s(n-1),s,\dots,s)}E_1)$ has the
property $t^{-s(n-1)}L_0 \supset L_{\wsn}\supset t^sL_0$ and
dim$\,_K(L/t^sL_0)=sn$ (here, $L_0$ is the standard $A$-lattice,
namely, the $A$-span of $\{e_1,\dots,e_n\}$). Further, for
$w=\tau^{(c_1,\dots,c_n)}\in \hW^{{\mathcal{P}}}$, it is easily
seen that $w\le\wsn$ if and only if the lattice
$L_{w}(=\tau^{(c_1,\dots,c_n)}E_1)$ has the property
$t^{-s(n-1)}L_0 \supset L\supset
t^sL_0\,\mathrm{and}\, \dim_K(L/t^sL_0)=sn$. Thus, we get an
identif\/ication: \[X(\wsn)=\{A{\texttt{-}}{\mathrm{lattice\,
}}L|t^{-s(n-1)}L_0 \supset L\supset
t^sL_0,\,\mathrm{and}\, \dim_K(L/t^sL_0)=sn\}.\]

 It is well known that $\calG/\calP = \varinjlim X(\wsn)$. We
therefore restrict our attention to the discussion of $X(\wsn)$.
Mainly for notational convenience we replace $X(\wsn)$ by
$t^{s(n-1)}X(\wsn) =: X(\ws)$ where $\ws =
\tau^{(0,sn,\dots,sn)}$.  Let $d = sn$ (which we will keep
throughout the text). As an element of $\Gr(\infty)$, $\wsn$
contains $E_{sn+1}$ and is contained in $E_{1-s(n-1)n)}$. Hence
$X(\wsn)$ embeds into $G(d,E_{1-d(n-1)}/E_{d+1})$. Consequently
$X(\ws)$ may be thought of as a subset of
\[t^{d-s}G(d,E_{1-d(n-1)}/E_{d + 1}) =
G(d,E_{1}/E_{dn+1}).\] We will denote $E_1/E_{dn +1}$ by $V_s$ or
simply $V$.

Let $u := 1 + t \in \GL(V)$, $u(v)=v+tv$, $v\in V$; then $u$ is
unipotent and clearly $X(\ws) \subseteq G(d,V)^u$. We shall now
show that this inclusion is in fact an equality.
\begin{prop}\label{kap0}
With notations as above, we have, $X(\ws) = G(d,V)^u$ \tu{(}the fixed
point set of $u$ with the reduced scheme structure\tu{)}.
\end{prop}
Before we can prove this proposition we need to introduce some notation also used throughout the rest of the paper.
Choosing the
basis on $V$ given by $e_1,e_2,\dots,e_{sn^2} \in E_1$, let $\TV
\subset \GL(V)$ be the induced diagonal torus and let $\BV$ be the
Borel subgroup of lower triangular matrices. Notice that $\calP$
acts on $V$ by means of a representation $\calP \to \GL(V)$, and
the image of $\calB$ is contained in~$\BV$. Similarly, $\hT$ acts
on $V$ and injects into $\TV$. By construction all the $\hT$-f\/ixed
points in $X(\ws)$ are actually $\TV$-f\/ixed points. The
$\TV$-f\/ixed points in $G(d,V)$ are the $d$-spans $\K e_{i_1} + K
e_{i_2} + \dots + K e_{i_d}$, where $1\leq i_1 < i_2 < \dots < i_d
\leq dn$. We shall denote such $d$-tuples by~$I_d$ or just~$I$. An
element $x = (x_1 < x_2 <\dots < x_d)\in I$ determines (uniquely)
a point of $\Gr(\infty)$, namely, the subspace of $K^{\infty}$
given by
\[
K e_{x_1} + \cdots + K
e_{x_d}+ \sum_{j>sn^2} K e_{j};
\] we shall denote it by
$E_{|x|}$ where
\[
\vert x \vert = x \cup \ZZ_{> sn^2}.
\] Note that $E_{|x|}$
contains $E_{sn^2 + 1}$. Now $E_{|x|}$ is in $\mathfrak X$  if and
only if the underlying space is $t$-stable if and only if $x$ is a
$u$-f\/ixed point if and only if for all $y\in |x|$, $y + n \in
\vert x\vert$ if and only if $x_i + n \in \vert x\vert$. Let $I^u$
denote the set of all $x$ with this property. We will identify an
element $x\in I^u$ with its counterpart $\tau^c \in \hW$.

Recall the Bruhat--Chevalley order $\succeq$ on $I$ with respect to $\BV$: it is def\/ined as $v \succeq w$ if $\overline{\BV v} \supseteq \BV w$. A similarly def\/ined order (now with respect to $\calB$) exists on the set of af\/f\/ine Schubert varieties.
For $I$ we also have the combinatorial partial order given by $v \geq w$ if and only if $v_i \leq w_i$  for all $i = 1,2,\dots,d$. A fundamental result in the theory of Schubert varieties is the fact that the two orderings on $I$ coincide, and moreover it is compatible with the order on the Weyl group of $\GL(V)$ given by the simple generators.

Therefore we obtain two orderings on $I^u$, one stemming from the Bruhat--Chevalley order on $\BV$-Schubert varieties, and one for the af\/f\/ine Schubert varieties parameterized by $I^u$. It is clear that the inherited order is a priori weaker than the Bruhat--Chevalley order (with respect to $\calB$) on $I^u$. However, if $v \geq w$ in $I^u$, then also $X(v) \supset X(w)$ which can be shown combinatorially, or geometrically by observing that for $v \in I^u$, $\calB v$ is dense in $\overline{\BV v}^u$.

\begin{proof}[Proof of Proposition~\ref{kap0}]
As mentioned above, the Bruhat--Chevalley partial order on $I^u$ can be
described combinatorially, and hence the partial order on Schubert varieties
 in $X(\wsn)$ can be described combinatorially. The $d$-tuple in
$I^u$ (namely, $(1,n+1,2n+1,\dots,(d-1)n+1)$) representing $\ws$
is the largest in $I^u$. Hence it follows that, as sets, $X(\ws) =
G(d,V)^u$.
\end{proof}

We will now simply write $\leq$ to denote the partial order on $I$ and $I^u$.
Note that as a consequence of the fact that the partial orders on $I^u$ coincide, for every $w \leq \ws \in I^u$, we have $X(w)
= Y(w)^u$, where for every $w \in I$, $Y(w)$ denotes the
$\BV$-Schubert variety $\overline{\BV w}$.

The Schubert varieties in $\calG/\calP$ share many properties with
their classical counterparts, both geometrically, and
combinatorially.

We f\/inish this section by a well known but nevertheless important lemma, whose proof in its current form was pointed out to us by one of the referees.

\begin{lem}\label{lem:linearizable}
The action of $\hT$ \tu{(}resp. $\TV$\tu{)} on $G(d,V)$ is locally
linearizable, that is, $G(d,V)$ is covered by open affine
$\hT$-stable \tu{(}resp. $\TV$-stable\tu{)} neighborhoods. Consequently the
same applies for every closed $\hT$-stable \tu{(}resp. $\TV$-stable\tu{)}
subvariety.
\end{lem}
\begin{proof}
The open $\BV$-orbit in $G(d,V)$ is both, af\/f\/ine and $\hT$- (resp. $\TV$-) stable.
Its $\WV$-translates cover all of $G(d,V)$.
\end{proof}

\section[Reflections and combinatorics]{Ref\/lections and combinatorics}\label{sec:REFLECTIONS}

Let
$d,V,Y(w)$ etc., be as in the previous section. Let $\WV$ be the
Weyl group of the pair $(\GL(V),\TV)$, and let $\mathcal{R}$
denote the set of ref\/lections in $\WV$. Let $x \in Y(w)$ be a
$\TV$-f\/ixed point. Set
\[
S_x(w)=\{r\in\mathcal{R}\,|\,x \neq rx \, (\text{in }G(d,V))
\text{ and }rx\le w\}.
\] Then we have (cf.~\cite{LSe}) that
\begin{gather}\label{eq*}
x
\mathrm{\ is\ a\ smooth\ point\ of\ } Y(w) \mathrm{\ if\ and\
only\ if\ } \#\,S_x(w) =\dim Y(w). 
\end{gather}  These ref\/lections
are in one-one correspondence with the $\TV$-stable curves in
$Y(w)$ containing $x$ (see \cite{CP}). In fact, if $r = r_\alpha$
is the ref\/lection associated to the root $\alpha \in \PhiV$, the
root system of $(\GL(V),\TV)$, and if $rx\neq x$, then either
$\overline{U_\alpha x}$ or $\overline{U_{-\alpha}x}$ is a
$\TV$-stable curve containing $x$ and $rx$. Here $U_\alpha$
denotes the one-dimensional unipotent group normalized by $\TV$
whose Lie algebra has $\TV$-weight $\alpha$. Using these results,
it is not hard to determine the singular locus of $Y(w)$.
\begin{defn}
The conjugate Meyer diagram $\CM(w)$ of $w$ (or $Y(w)$) is the
Young diagram with $d$ rows whose $i$-th row consists of $d(n-1) +
i -w_i$ boxes.
\end{defn}
Notice that $d(n-1) + i$ is the $i$-th entry of the unique
$\BV$-f\/ixed point $e \in G(d,V)$. Clearly $e \in I^u$ is also the
unique $\calB$-f\/ixed point in $X(\ws)$. The maximal singularities
(maximal with respect to the Bruhat--Chevalley order) then are
given by the \emph{hooks} of $\CM(w)$ as follows. A~hook~$\calH$
is a~sequence of consecutive rows $R_i,R_{i+1},\dots,R_{i+k}$ of
$\CM(w)$ ($k > 0$) such that for the length $\vert R_j \vert$ of
row $R_j$ we have
\[ \vert R_i\vert > \vert R_{i+1}\vert = \vert R_{i+2}\vert = \dots = \vert
R_{i+k}\vert > \vert R_{i+k+1}\vert,
\] where by convention
$R_{d+1}$ is an empty row in case $i + k = d$. The element
$w_\calH \leq w$ is the unique element of $I$, such that
$\CM(w_\calH)$ is obtained from $\CM(w)$ by replacing the rows
$R_i,R_{i+1},\dots,R_{i+k}$ by rows of equal length $\vert R_{i+1}
\vert -1$. Equivalently, $w_\calH$ is obtained from $w$ by
replacing $w_i$ with $w_{i+k} + 1$. It is well known that $w_\calH$ is a singularity of $Y(w)$ (see for instance \cite{LW} for a more general result).

In the af\/f\/ine setting however, things are more complicated.
Firstly, \eqref{eq*} need not hold. Se\-condly, it is more complicated to
even describe the smooth $\hT$-f\/ixed points in combinatorial terms.

\begin{defn}\label{si} Let $x \in I^u$. For $1\leq i \leq n$,
we def\/ine \[ S_i(x) = \{ j \in x \mid j \equiv i \mod n\}\] the
$n$-\emph{string} through $i$ in $x$. \end{defn}

Denote
$\ell_i(x): = \vert S_i(x)\vert$. Notice that if $x = \tau^c$ and
$c = (c_1,c_2,\dots,c_n)$, then $\ell_i(x) = d-c_i$. Def\/ine
$h_i(x)$ as the minimal element of $S_i(x)$, if $S_i(x)$ is
nonempty (thus, $h_i(x)$ is the ``head'' of the string), and $h_i(x)
= dn + i$, if $S_i(x)$ is empty. (Note that $h_i(x) = i + c_in$).
It will be convenient to denote the unique integer between $1$ and
$n$ congruent to a given integer $i$ mod $n$ by $[i]$; thus, any
$x_k \in x$ satisf\/ies $x_k \in S_{[x_k]}(x)$.

Let $\halpha \in \hPhi$. Then clearly $s_\halpha(\vert x\vert)$ is
again a sequence of integers.
\begin{defn} We say $s_\halpha$ is \emph{defined
at $x$}, if $s_\halpha (\vert x \vert) = \vert y \vert$ for some
$y \in I^u$. \end{defn} Note that this is equivalent to saying
that $s_\halpha x \leq \ws$. If $s_\halpha$ is def\/ined at $x$,
then it operates on the strings, that is $s_\alpha x$ is obtained
from $x$ by removing a number of elements in $S_i(x)$ and adding
the same number to another string $S_j(x)$. Thus, $s_\halpha x$ is
determined by the requirements that $\ell_i(s_\alpha x) =
\ell_i(x) - k$ and $\ell_j(s_\halpha x) = \ell_j(x) + k$ for some
suitable $k \geq 0$. The indices $i$ and $j$ are referred to as
the \emph{indices corresponding to $s$}. The simplest way to
describe these operations is by means of the following diagram:
\begin{defn}
Let $x \in I^u$. The \emph{string diagram} $\Sigma(x)$ of $x$
consists of $n$ rows where the $i$-th row has $\ell_i(x)$ boxes.
\end{defn}
Obviously $\Sigma(x)$ is just an encoding of $c \in \ZZ^n$ in the
description of $x$ as $x = \tau^c$, and $\ell_i(x) = d-c_i$. The
total number of boxes in $\Sigma(x)$ is always $d$.

\begin{figure}[t]
\centerline{\includegraphics{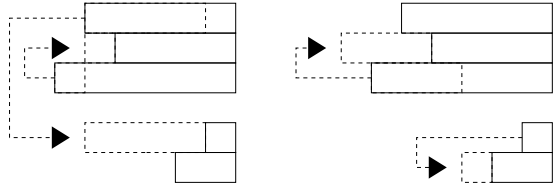}}
\caption{An example of $\Sigma(x)$ where $n
= 6$, $s = 3$, $d=18$, $x =
(\mathbf{75},\mathbf{79},81,85,\mathbf{86},87,91, 92$, $93,97,98,99,\mathbf{102},103,104,105,107,\mathbf{108})$
(the heads of the strings are in bold face), together with two
small down-exchanges (left) and two large up-exchanges (right).}\label{fig:REFLECTIONS}
\end{figure}

\subsection[Small reflections]{Small ref\/lections}

We shall denote a typical ref\/lection
in $\hW$ by $s$, though $s$ has also been used as a superscript in
$\ws$. But, we believe, this will not create any mix-up; whenever
necessary, we will be explicit about the particular reference.
There is a special class of ref\/lections, which will play a crucial
role in our description of singularities.
\begin{defn}
Let $x \in I^u$, and let $s\in \hW$ be a ref\/lection def\/ined at
$x$, such that $sx \neq x$.  Then~$s$ is called \emph{small}, if
and only if $\vert s(x_k) - x_k\vert < n$ for all $1 \leq k \leq
d$. Otherwise, $s$ is called \emph{large}.
\end{defn}
\begin{rem}
Let $s = s_\halpha$ be any ref\/lection def\/ined at $x$ with
$s_\halpha x \neq x$. If $s_\halpha x > x$ (respectively
$s_\halpha x < x$), there is a unique small ref\/lection $s' =
s_{\halpha'}$ with $x < s'x \leq sx$ (respectively, $x > s'x \geq
sx$), and $\Re(\halpha) = \Re(\halpha')$ (cf.\ Section~\ref{weyl}).
As an example, we treat the case $s_{\halpha} x > x$. Suppose $\Re(\halpha) = (ij)$ with $i < j$, and given $\ell_i(x) \geq \ell_j(x)$, then $\halpha' = (ij) + \delta$. If $\ell_i(x) < \ell_j(x)$, then $\halpha' = (ij)$.
\end{rem}

The small ref\/lections are easily described in terms of $\Sigma(x)$
(see Fig.~\ref{fig:REFLECTIONS}). For any pair of integers $1
\leq i < j \leq n$ there is almost always a unique small
ref\/lection $s$ with $sx > x$ and corresponding indices $i$ and
$j$: $sx$ is obtained from $x$ by $\ell_i(sx) = \ell_j(x)$ and
$\ell_j(sx) = \ell_i(x)$ if $\ell_j(x) > \ell_i(x)$ (i.e.\ the
rows of $\Sigma(x)$ at positions $i$ and $j$ are simply switched);
if on the other hand $\ell_i(x) \geq \ell_j(x)
> 0$, then $sx$ satisf\/ies $\ell_i(sx) = \ell_j(sx)-1$, and $\ell_j(sx) =
\ell_i(x) + 1$. The only case when $s$ does not exist is
$\ell_j(x) = 0$. We refer to the process of applying $s$ as
\emph{up-exchanging~$i$ and~$j$}.

Similarly, \emph{down-exchanging $i$ and $j$} is the inverse
procedure, i.e.\ the result of down-exchanging $i$ and $j$ is the
unique $x' < x$, such that up-exchanging $i$ and $j$ in $x'$ gives
$x$. If $\ell_i(x) = \ell_j(x)$ or if $\ell_i(x) = \ell_j(x)-1$,
then $x'$ is not def\/ined in this manner, and we let $x' = x$, in
this case.

To simplify our notation and get rid of the two dif\/ferent cases
when up-exchanging (or down-exchanging) we make a def\/inition:
\begin{defn}\label{lx} Def\/ine
\[L(x) = (\ell_1(x),\ell_2(x),\dots,\ell_n(x),\ell_1(x) +
1,\ell_2(x) + 1,\dots, \ell_n(x) + 1).
\] \end{defn} Notice that $x$
is uniquely determined by any $n$ consecutive entries of $L(x)$
(together with the f\/irst entry).

Let $L(x) = (l_1,\dots,l_{2n})$. Up/down-exchanging $i$ and $j$
for $1\leq i,j\leq 2n$ always refers to up/down-exchanging $[i]$
and $[j]$. Suppose we want to up-exchange $i < j \leq n$. By
replacing $i$ with $i + n$ and switching $i$ and $j$ if necessary,
we may assume that $i < j$ and $l_i < l_j$ with $i \leq n$. If we
now def\/ine~$x$ by the~$n$ entries $l_i,\dots,l_{i+n-1}$, then the
corresponding entries in $L(sx)$ are the same with the exception
that $l_i$ and $l_j$ switch positions. We will often describe $x$
by any~$n$ entries in $L(x)$ which contain $l_i$ or $l_{i+n}$ for
every $1 \leq i \leq n$. Similarly, down-exchanging~$[i]$ and~$[j]$ reduces to switching the positions of $l_i$ and $l_j$ if we
pick~$i$ and~$j$ such that~$i < j$ and~$l_i > l_j$. We write $i
\da j$ (resp. $i \up j$) for the down-exchange (resp. up-exchange)
of~$i$ and~$j$.

\subsection[The codimension of $X(sx)$ in $X(x)$, $sx<x$]{The codimension of $\boldsymbol{X(sx)}$ in $\boldsymbol{X(x)}$, $\boldsymbol{sx<x}$}

Recall that a property of the Bruhat--Chevalley ordering is the
fact that for any $x \leq w$ we have $\codim_{X(w)}(X(x)) = \max
\{ i \mid \exists \, x = \tau_0 < \tau_1 < \dots < \tau_i = w \}$.
This still holds in the af\/f\/ine Grassmannian (easily proven
combinatorially, or by the fact that the $\calB$-orbits are
isomorphic to af\/f\/ine spaces).

One reason for introducing $L(x)$ is the fact that if $s$ is a
small ref\/lection and $sx < x$, the codimension of $X(sx)$ in
$X(x)$ may be read of\/f immediately, thanks to the following lemma:
\begin{lem}\label{lem:CODIM}
Let $x \in I^u$ and let $s$ be a small reflection with $sx < x$.
If $sx$ is the result of down-exchanging $i < j < i + n$ with $l_i
> l_j$ and $i < n$, say, then the codimension of $X(sx)$ in~$X(x)$
equals $1 + g_\geq + g_>$ where
\[ g_\geq = g_\geq(i,j,x) = \vert \{ i < k < j\mid L(x)_i \geq L(x)_k \geq L(x)_j\} \vert \]
and
\[ g_> = g_>(i,j,x) = \vert \{ i < k < j \mid L(x)_i > L(x)_k > L(x)_j\}.\]
\end{lem}
Notice that the assumptions $i < j$ and $L(x)_i > L(x)_j$ are no
restriction as we may replace $i$ or $j$ if necessary with $i+n$
or $j+n$, respectively.
\begin{proof}
We proceed by induction on $g_\geq$. First suppose $g_\geq = 0$.
Then we have to show that $sx$ has codimension one. Equivalently,
if $sx \leq y \leq x$, then either $y = x$ or $y = sx$. Suppose
such a~$y$ is given. We may describe $x$, $y$, and $sx$ by the
entries of $L(\cdot)_{i+k}$ for $k = 0,1,\dots,n-1$.

First observe that $x$ and $sx$ coincide before (and including)
the entry preceding $h_i(x)$ in $x$ and after (and including) the
position of $h_i(sx)$ in $sx$. This immediately shows that $L(x)_i
\geq L(y)_i \geq L(sx)_i$. Similarly, $x$ and $sx$ coincide after
the position of $h_j(x) - n$ in $sx$. Thus $L(y)_j \geq L(x)_j$.
Now, if $L(y)_j > L(sx)_j = L(x)_i$, then $h_j(y) < h_i(x)$ and so
would be in a range where $x$ and $sx$ coincide, a contradiction.
Thus $L(sx)_j \geq L(y)_j \geq L(x)_j$.

For $k \neq i,j$, it follows that $L(y)_k = L(x)_k = L(sx)_k$: indeed, suppose f\/irst that $i < k < j$; then as $h_{[k]}(x) < h_i(x)$ or $h_{[k]}(x) > h_j(x)$ for such $k$ (in view of the hypothesis that $g_{\geq} = 0$), the position of $h_{[k]}(x)$ is in a range where $x$ and $sx$ agree. Moreover, $h_{[x]}(x) - n$ cannot be present in $y$, as it is not present in either of $x$ and $sx$, and would also fall into a range where $x$ and $sx$ agree (being strictly bigger than $h_j(x) - n$). Consequently, $h_{[k]}(y) = h_{[k]}(x)$, and it must have the same position.

Similarly, if $k > j$, and suppose $y$ contains $h \equiv k
\mod n$. Pick $r < h$ maximal congruent $i$. Then there is $t
\equiv  j \mod n$ such that $r < t < h < r + n$. There are two
cases: either $r$ is present in $x$, and $t$ is not, in which case
$r$ is replaced by $t$ in $sx$, and the entries in $x$ and $sx$
strictly between $t$ and $r+n$ coincide and don't change
positions; or, both $r$ and $t$, or neither of them, is present
in $x$. In which case $x$ and $sx$ coincide at values strictly
between $r$ and $r + n$ (including positions). It thus follows
that $h$ must be in $x$ and $sx$. Consequently $h_{[k]}(y)
= h_{[k]}(x)$ for all these $k$ (at the same position) and
therefore $L(y)_{k} = L_{k}(x) = L_{k}(sx)$.

Since the coordinate sums of $L(y)$, $L(x)$, and $L(sx)$ all
agree, it follows that $L_i(y) + L_j(y) = L_i(sx) + L_j(sx) =
L_i(x) + L_j(x)$. If $L_i(x) > L_i(y) > L_i(sx)$,  then
necessarily $L_j(y) < L_j(sx)$, and the number of elements in $y$
that are less than or equal to $h_j(sx)$ is strictly smaller than
the same number for $sx$, contradicting $y \geq sx$. The only two
possibilities now are $L_i(y) = L_i(sx)$ or $L_i(y) = L_i(x)$,
resulting in $y = sx$ or $y=x$ as claimed. This completes the
proof in the case $g_{\geq} = 0$.

\begin{figure}[t]
\centerline{\includegraphics{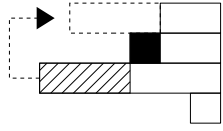}}
\caption{Here $n = 4$, and $1$ and $3$
(resp.~$3$ and $5$) are down-exchanged. The black box marks the
entry in~$x$ and $sx$ from which on they coincide. Since $g_{\geq}
= 0$, this is a codimension one down-exchange.}\label{fig:CODIM_POSITION}
\end{figure}

Let then $g_{\geq} > 0$, and let $i'$ be the f\/irst index greater
than $i$ such that $l_i\geq l_{i'} \geq l_j$, then clearly $x \geq
y \geq sx$, where $y$ is the result of $i \da i'$  (if $l_i =
l_{i'}$, then $x = y$): $sx$ is obtained from $y$ by $i' \da j$,
then $i \da i'$. Notice that the last down-exchange is necessary
only if $L(x)_{i'}>L(x)_j$. Consider the down-exchange of $i'$ and
$j$ in $y$, resulting in $y' \geq sx$. Clearly $g_{\geq}(i',j,y) =
g_\geq(i,j,x)-1$ ($L(y)_{i'} = L(x)_i$). By induction, the
codimension of $y'$ in $X(y)$ is $1 + g_\geq(i',j,y) +
g_>(i',j,y)$. In addition, $g_>(i',j,y) = g_>(i,j,x)$ if $y = x$,
and $g_>(i',j,y) = g_>(i,j,x)-1$ otherwise. If $y' \neq sx$, then
the codimension of $sx$ in $X(y')$ is $1$ by construction and the
application of the case $g_\geq = 0$. Similarly, if $y \neq x$
then the codimension of $y$ in $X(x)$ is one as well. The result
now follows.

One caveat is the following subtlety: strictly speaking, applying
the induction hypothesis requires $i' \leq n$. However, if $i' >
n$ the numbers $g_>(i',j,y)$ and $g_\geq(i',j,y)$ do not change if~$i'$ and~$j$ are replaced with $i'-n$ and $j-n$.
\end{proof}

\begin{figure}[t]
\centerline{\includegraphics{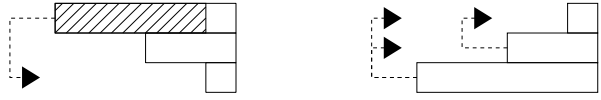}}
\caption{A down-exchange of
$1$ and $3$ (left) and the corresponding ref\/lections (right).
Notice that here~$g_\geq$ and~$g_>$ both are $1$.}\label{fig:REFLECTIONS_DOWNEXCHANGE}
\end{figure}

\begin{rem}\label{rem:CORRESPONDING_REFLECTIONS}
Suppose $x$, $s$, $i$ and $j$ are as in Lemma~\ref{lem:CODIM}.
Notice that an element of $G_\geq(i,j,x) := \{ i< k < j\mid L(x)_i
\geq L(x)_k \geq L(x)_j\}$ gives rise to either one or two
ref\/lections $s'$ with $x > s'sx > sx$ as follows: If $L(x)_i
> L(x)_k$ up-exchanging $j$ and $k$ in $sx$ is possible and the
result is below $x$, because it is obtained from $x$ by $k \da j$
and then $i \da k$. If on the other hand $L(x)_k > L(x)_j$, then
up-exchanging $k$ and $i$ is possible in $sx$ and below $x$,
because it is the same as  $k \da j$ and $i \da j$ applied to $x$.

Thus there are precisely $g_>$ elements of $G_\geq(i,j,x)$ giving
rise to two ref\/lections. The total number of ref\/lections thus
obtained including $s$ is $1 + g_\geq + g_>$, the codimension of
$sx$ in $x$. We refer to these ref\/lections as the ref\/lections
\emph{corresponding to the down-exchange} of $i$ and $j$. For an
example see Fig.~\ref{fig:REFLECTIONS_DOWNEXCHANGE}.
\end{rem}

\section[The connection between $X(w)$ and $Y(w)$]{The connection between $\boldsymbol{X(w)}$ and $\boldsymbol{Y(w)}$}\label{sec:RELATION}

In this section we will further investigate the connection between
$X(w) \subset Y(w)$. One might think that $X(w) = Y(w)^u$ also as
a scheme, but as it turns out this is not true in general, as it
may happen that $T_x(Y(w))^u \supsetneq T_x(X(w))$ at some $x \leq
w$. Nevertheless, a f\/irst step in computing $T_x(X(w))$ is to
determine $T_x(Y(w))^u$, and in some cases knowledge of the latter
is enough to determine the former.

\subsection{Tangents to the Grassmannian}
We will need the following well
known description of the tangent space to $T_x(G(d,V))$ at a~point~$x$. Viewing $x$ as a~subspace of $V$, let  $p_x\colon V \to V/x$
and $i_x \colon x \to V$ be the projection and inclusion maps,
respectively.
\begin{lem}\label{lem:GRASS_TANGENT}
Let $x\in G(d,V)$. Then \[ T_x(G(d,V)) = \Hom(x,V/x)\] in a
natural way. In addition, the differential of the orbit map
$\GL(V) \to G(d,V)$ which sends $g$ to $gx$, is given by
\[
\begin{array}{ccc}
\End(V) & \to & \Hom(x,V/x),\\
\xi & \mapsto & p_x \xi i_x.
\end{array}
\]
If $P$ denotes the stabilizer of $x$ in $\GL(V)$, this map is
equivariant with respect to the adjoint action of $P$ on
$\End(V)$, and the natural action of $P$ on $\Hom(x,V/x)$.
\end{lem}
\begin{rem}\label{wt} If $E_{ij}$ ($1\leq i,j \leq dn$) denotes the element
in End$\,V$, sending $e_j$ to $e_i$ and if $x =
(x_1,x_2,\dots,x_d) \in I$ , then $T_x(G(d,V))$ is spanned by the
images of those $E_{ij}$ for which $j \in x$ but $i \not \in x$.
We will denote these elements of $T_x(G(d,V)$ by $E_{ij}$ as well.
Notice that $E_{ij}$ is a $\TV$-eigenvector of weight $\epsilon_i
- \epsilon_j$ (where $\epsilon_k$ is the element in the character
group of $\TV$, sending a~diagonal matrix in $\TV$ to its $k$-th
diagonal entry). Thus, its $\hT$-weight is $\epsilon_{[i]} -
\epsilon_{[j]} + h\delta$ for a~suitable~$h$, where either
$[i]\neq [j]$ or $h \neq 0$. In particular, all the $\hT$-weights
of $T_x(G(d,V))$ are roots.
\end{rem}

\subsection{Real tangents}
It is now clear, that if $x$ is an element of $I^u$, then
$\xi \in T_x(G(d,V))$ is $u$-f\/ixed if and only if $\tau \xi = \xi
\tau$. From this it follows easily that any $\xi \in
T_x(G(d,V))^u$ is uniquely determined by its values on~$e_{h_i(x)}$ ($1 \leq i \leq n $) such that $h_i(x) \leq dn$; for,
then $\xi(e_{h_i(x) + rn}) = \tau^r\xi(e_{h_i(x)})$, $0 \leq r <
\ell_i(x)$.

\begin{defn}\label{def:REAL_TANGENT}
Let $x \in I^u$, and suppose $\halpha = \alpha + h\delta \in
\hPhi$. We may write $\alpha = (ij)$ for some $1 \leq i,j \leq n$.
If $0 < \ell_j(x)$ and $h_j(x) = j + tn$ such that $t + h \geq 0$,
$i + (t+h)n < h_i(x)$, and $i + (t + h + \ell_j(x) )n  \geq
h_i(x)$, then $\xi_\halpha$ in $\Hom(x,V/x)$ is def\/ined as:
\begin{equation*}
\xi_\halpha e_{k} = \begin{cases} p_x e_{i + (r + t + h)n}, \quad
&r > 0,\  k = h_j(x) + rn .\\
0, & \text{otherwise}.
\end{cases}
\end{equation*}
To avoid having to state the hypotheses over and over again, we
simply say, \emph{$\xi_\halpha$ is defined at $x$} to indicate
that all conditions above are met.
\end{defn}
\begin{lem}Suppose $\xi_\halpha$ is defined at $x \in I^u$. Then
$\xi_\alpha$ is a $u$-fixed $\hT$-eigenvector of $T_x(G(d,V))$ of
weight $\halpha$. Conversely, if $\xi$ is any $\hT$-eigenvector in
$T_x(G(d,V))^u$ whose weight is a real root~$\halpha$, then
$\xi_\halpha$ is defined at $x$ and $\xi \in \K \xi_\halpha$.
\end{lem}
\begin{proof}
The f\/irst statement is immediate from the def\/inition; we will
therefore prove the second assertion: Suppose $\xi \in
T_x(G(d,V))^u$ has $\hT$-weight $\halpha = \alpha + h\delta$.
Again write $\alpha = (ij)$. Let $x_r \in S_k(x)$ (cf. Def\/inition
\ref{si}), $x_r = k + mn$, say. Then the $\hT$-weight of $e_{x_r}$
is $\epsilon_k + m \delta$. As $\xi e_{x_r}$ has weight
$\epsilon_k + \alpha + (m+ h)\delta$, and by the description of
possible $\hT$-weights in $\Hom(x,V/x)$ (cf. Remark \ref{wt}) it
follows that $k$ must equal $j$, if $\xi e_{x_r}$ is not to be
zero. Therefore $\xi$ is uniquely determined by its values at
$e_{h_j(x)}$, and it is suf\/f\/icient to show that $\xi_\halpha$ is
def\/ined at $x$ and that $\xi e_{x_r} \in K \xi_\halpha e_{x_r} = K
e_{i + (h + m)n}$ in case $x_r = k + mn = h_j(x)$. But $e_{i+
(h+m)n}$ is the only possibility of an element of $V/x$ with
weight $\epsilon_i + (h+m)\delta$ (up to scalars). As $\xi$ is
nonzero, it follows $\xi e_{h_j(x)} = c e_{i + (h+m)n}$ for some
$c \neq 0$, and in particular $0 \leq (h + m)$ and $i + (h+m)n <
h_i(x)$. Since $\xi$ is also $\tau$-equivariant we must have
$\tau^{\ell_j(x)}\xi e_{h_j(x)} = 0$, and thus $i +
(h+m+\ell_j(x))n \geq h_i(x)$. But now $\xi_\halpha$ is def\/ined at
$x$ and obviously $\xi$ is proportional to $\xi_\halpha$.
\end{proof}

We will refer to the elements of $T_x(G(d,V))^u$ which
have a real root as $\hT$-weight as \emph{real tangents}.  Notice that
if $\xi_\halpha$ is def\/ined at $x$, and $\halpha > 0$, then
$\xi_\halpha$ actually lifts to a $\tau$-invariant element of
$\End(V)$ of $\hT$-weight $\halpha$. In fact, if $\halpha = (ij) +
h\delta > 0$ (and in particular, $h\geq 0$) then
\begin{equation} \label{eq:END_XI}\xi_\halpha = E_{i +
hn,j} + E_{i+(h+1)n,j+n}+ \dots + E_{i+(d-1)n,j+(d-h-1)n} \in
\End(V)^u.
\end{equation}
If $h < 0$ this is not possible (as then $j + (d-h-1)n > dn$, so
no element of $V$ can be mapped to $e_{i + (d-1)n}$). But we still
have
\begin{equation}\label{eq:XI_ALPHA}
\xi_\halpha = \sum_{\substack{0 \leq r < d\\ j + rn \geq h_j(x) }}
E_{i+ (r+h)n, j+ rn}
\end{equation}
in $\Hom(x,V/x)$ with the convention that $E_{i,j} = 0$ if $i$ or
$j > dn$.

\subsection[Reflections and $\hT$-curves]{Ref\/lections and $\boldsymbol{\hT}$-curves}

The main goal of this section is to show that actually $\xi_\halpha \in T_x(X(\ws))$ whenever $\xi_\halpha$ is a real tangent def\/ined at $x$.
As noted in the previous section, the $\TV$-stable curves play a
crucial role when determining the singularities of a classical
Schubert variety. They  still give rise to a~ne\-ces\-sary though not
suf\/f\/icient criterion in the case of af\/f\/ine Schubert varieties. Let
\mbox{$x\leq w\in I^u$}. A~\mbox{\emph{$\hT$-curve}} through $x$ in $X(w)$ is
the closure of a one-dimensional $\hT$-orbit in $X(w)$ which
contains~$x$. We denote the set of all $\hT$-curves through $x$ by
$E(X(w),x)$. By results of~\cite{CP}, each such $\hT$-curve in $E(X(w),x)$ is the
$G_\halpha$-orbit of $x$ for some suitable $\halpha \in \hPhi$,
where $G_\halpha$ is the copy of $\SL_2(\K)$ in $\calG$ which is
generated by the root-groups $U_{\pm \halpha}$; here, for any
$\halpha \in \hPhi$, $U_\halpha$ is the (uniquely determined)
image of an inclusion $x_\halpha \colon \K \to \calG$ that is
equivariant with respect to the $\hT$-actions ($\hT$ acts on
$\K$ by $\halpha$ and by conjugation on~$\calG$):
$x_\halpha(\halpha(t)k) = t x_\halpha(k) t^{-1}$ (conjugation here
means that $t = (t_0,s) \in \hT = T \times K^*$ acts as $tgt^{-1}
:= st_0gt_0^{-1}$)). It follows that around $x$, any such~$C$ has
the form $U_\halpha x$ for a suitable $\halpha$, and $T_x(C)$ is a
line in $T_x(X(w))$ with $\hT$-weight~$\halpha$. In particular $C$
is smooth. Moreover only real roots occur as weights. Let
$TE(X(w),x))$ denote the ``span'' of the $\hT$-curves, that is,
\[TE(X(w),x) = \bigoplus_{C\in E(X(w),x)} T_x(C).\] An immediate
consequence of Lemma~\ref{lem:linearizable} is the fact that
$\vert E(X(w),x) \vert \geq \dim X(w)$ (for a proof see~\cite{CP}). This is
sometimes referred to as Deodhar's inequality. As no two
$\hT$-curves have the same $\hT$-weight at $x$, it follows that
$\dim TE(X(w),x) \geq \dim X(w)$. Summarizing, let us recall the following necessary criterion from \cite{CP} for $x \leq w$ being a smooth point of $X(w)$:
\begin{lem}\label{lem:TOOMANYCURVES}
Let $x \leq w$; if $x$ is a smooth point of $X(w)$, then
\[ \vert E(X(w),x) \vert = \dim X(w). \]
Equivalently, the number of reflections $s$ such that $sx \neq x$
and $sx \leq w$ equals $\dim X(w)$.  This in turn is equivalent to
\[ \vert \{ s \mid x < sx \leq w \} \vert = \codim_{X(w)}(X(x)). \]
\end{lem}
The last statement of the Lemma is an immediate consequence of the
fact that $\vert E(X(x),x) \vert = \dim X(x)$, as $x$ is a smooth
point of $X(x)$.

It should be pointed out, however, that contrary to the classical setting, this condition is not suf\/f\/icient (see Remark~\ref{rem:Peterson}).

One of the reasons that smoothness of Schubert varieties in the
af\/f\/ine Grassmannian is a~more delicate question than in the
ordinary Grassmannian, is the existence of imaginary roots. For
instance we have:
\begin{lem}\label{lem:IMAGINARY_TANGENT}
Let $x \leq w \in I^u$. If there is a line $L \subset T_x(X(w))$
which is a $\hT$-eigenvector whose weight is an imaginary root,
then $x$ is a singular point.
\end{lem}
\begin{proof}By the remarks preceding the lemma, the $\hT$-weights
of $TE(X(w),x)$ are real roots. Thus, $L \neq T_x(C)$ for all $C
\in E(X(w),x)$. The lemma now follows from the following
well-known fact: If a  torus $S$ acts on an af\/f\/ine variety $X$
with smooth f\/ixed point $x$, then for every $S$-stable subspace $M
\subset T_x(X)$ there exists an $S$-stable subvariety $X' \subset X$ such that
$T_x(X') = M$. In the case of $X = X(w)$, applying this to an open
af\/f\/ine neighborhood of $x$, and putting $M = L$, the result
follows.
\end{proof}
\begin{rem}\label{preceed} While it is true that $\xi_\halpha$ is actually tangent to the
$\hT$-curve $\overline{U_\halpha x} \subset \calG/\calP$ it is not
always possible to realize this identif\/ication inside $G(d,V)$,
since $U_\halpha$ may not act on $V$ (in particular if $\halpha <
0$) commuting with $\tau$. If $\halpha$ is positive, then
$\xi_\halpha \in \End(V)$ is nilpotent and actually spans the
image of $\Lie(U_\halpha)$ in $\End(V)$. Consequently,
$U_\halpha\subset \calB$ injects into $\GL(V)$. In fact, let
$U_{kl} \subset \GL(V)$ denote the root group with Lie algebra $\K
E_{kl}$. Then the image of $U_\halpha$ is a~one-dimensional
subgroup of
\[U:= U_{i+hn,j}U_{i+(h+1)n,j+n}\cdots U_{i+(d-1)n,j+(d-h-1)n}.\]
Notice that all the individual factors in this product mutually
commute, and that this product therefore is direct and a subgroup
of $\GL(V)$, and $U_\halpha = U \cap C(u)$ where $C(u)$ denotes
the centralizer of $u$ in $\GL(V)$ (here, $u=1+\tau$ (cf. Section~\ref{iden})).
\end{rem}
\begin{lem}\label{lem:REFLECTION_TANGENT}Let $\halpha \in \hPhi$. If $\halpha > 0$, then~$\xi_\halpha$ is defined at $x$ if and only if $s_\halpha x < x$.
If $\halpha < 0$, then~$\xi_\halpha$ is defined at $x$ if and only
if $x < s_\halpha x \leq \ws$.

As a consequence, $\xi_\halpha \in TE(X(\ws),x)$, if it is
defined.
\end{lem}
\begin{proof}
If $\halpha > 0$, and $\xi_\halpha$ is def\/ined at $x$ then by the
remarks preceding the lemma, $U_\halpha \subset \GL(V)$, and
$U_\halpha x \neq x$. Since $\overline{U_\halpha x}$ is a
$\hT$-curve connecting $x$ and $s_\halpha x$, the result follows.

If $\halpha < 0$, then the conditions for $\xi_\halpha$ to be
def\/ined at $x$ assert that $x < s_\halpha x \leq \ws$. Then $C :=
\overline{U_{-\halpha}s_\halpha x}$ is a $\hT$-curve of $X(\ws)$
containing $x$. Since $\ker \halpha \subset \hT$
acts trivially on this curve, its tangent lines must have
$\hT$-weights in $\QQ \halpha \subset X(\hT)\otimes \QQ$.
Obviously, $\halpha$ is the only $\hT$-weight of $T_x(G(d,V))^u$
satisfying this condition, and the only corresponding eigenvector
is $\xi_\halpha$. As a~consequence $T_x(C) = \K \xi_\halpha
\subset TE(X(\ws),x)$.
\end{proof}
If for any $\hT$-stable subspace $M \subset T_x(G(d,V))$, $M_{\real}$
denotes the span of weight-subspaces for real roots, then we have
seen:
\begin{cor}\label{cor:CURVE_COR}
For any $x \leq \ws$ we have \[ T_x(X(\ws))_{\real} = TE(X(\ws),x)
= T_x(G(d,V))^u_\real.\]
\end{cor}
For general $w$ the situation is more delicate. One might think
that $T_x(X(w))_{\real} = T_x(Y(w))^u_\real$, but this is not true
in general. Also it is not clear whether $TE(X(w),x) =
T_x(X(w))_\real$.

\begin{rem}\label{rem:REALS_IN_Y}
Let $x, w\in I^u,x \leq w$. Suppose that $\xi_\halpha$ is def\/ined
at $x$ for $\halpha = (ij) + h\delta$. Recalling (cf.\
equation~\eqref{eq:XI_ALPHA}) that $\xi_{\halpha} \in T_x(Y(w))^u$
if and only if $E_{i+(k+h)n,j+kn} \in T_x(Y(w))$ for all~$k$ such
that $0 < i + (k+h)n < h_i(x)$ and $k < d$. Indeed, as $T_x(Y(w))$
is $\TV$-stable, $\xi_\halpha$ is contained in $T_x(Y(w))$ if and
only if every $\TV$-eigenvector it is supported in is an element
of $T_x(Y(w))$. This in turn is equivalent to saying that
$r_{i+(k+h)n,j+kn} x \leq w$ for all such $k$, where $r_{pq}$
denotes the transposition $(pq)$ in $\Sc_{dn} = \WV$. Notice that
all these $r_{i+(k+h)n,j+kn}$ commute.
\end{rem}

Let $k_0$ be the maximal $k$ appearing in
equation~\eqref{eq:XI_ALPHA}. Then $s_\halpha x \leq w$ if and only
if
\[r_{i+hn,j}r_{i+(1+h)n,j+kn}\cdots r_{i+(k_0 + h)n,j+k_0n} x \leq
w\] a stronger condition than having just $r_{i+(k+h)n,j+kn}x \leq
w$ for all $0 \leq k \leq k_0$ (at least for negative~$\halpha$).
If $s_\halpha$ is small however, the situation is dif\/ferent.
Keeping, the notation just introduced, we have:

\begin{lem}\label{lem:SMALL_REFLECTION}Let $s_\halpha$ be a small reflection $($such that $\xi_\halpha$ or $\xi_{-\halpha}$ is defined at $x \leq w)$. The following are equivalent:
\begin{enumerate}\itemsep=0pt
\item[{\rm 1)}] $s_\halpha x \leq w$;
\item[{\rm 2)}] $r_{i+{h+k}n,j+kn}x \leq w$ for all $0 \leq k \leq k_0$;
\item[{\rm 3)}] $\xi_\halpha$ or $\xi_{-\halpha}$ is tangent to
$X(w)$;
\item[{\rm 4)}] $\xi_\halpha$ or $\xi_{-\halpha}$ is tangent to $Y(w)$.
\end{enumerate}
\end{lem}
\begin{proof}
There is nothing to show if $s_\halpha x < x$. So assume
$s_\halpha x > x $ and $\halpha < 0$. If $s_\halpha$ is small,
then the positions where $r_{i+(h+k)n,j+kn} x$ and
$r_{i+(h+k'),j+k'n}x$ dif\/fer from $x$, are disjoint intervals in
$[1,d]$. Consequently, $s_\halpha x \leq w$ if and only if
$r_{i+(h+k),j+kn}x \leq w$ for all $k = 0,1,\dots,k_0$. The latter
condition is in turn equivalent to the fact that $\xi_\halpha \in
T_x(Y(w))$. Since $\xi_\halpha \in T_x(X(w)$ always implies
$\xi_\halpha \in T_x(Y(w))$ the lemma now follows.
\end{proof}

\subsection{Imaginary tangents}
Consistent with the notation introduced above, we call a tangent
$\xi \in T_x(G(d,V))^u$ \emph{imaginary}, if it is a
$\hT$-eigenvector for an imaginary root.
The Weyl group $\hW$ f\/ixes the imaginary roots $\ZZ\delta$ identically, i.e.\ $w(\delta) = \delta$ for all $w \in \hW$. Since for $x = e\calP$, the set of $\hT$-weights of $T_x(\calG/\calP)$ does not contain any positive imaginary roots, this means the same applies at any $\hT$ f\/ixed point $x \in \hW e\calP$, as the weights at $x$ are the $\hW$-translates of the weights at $e\calP$.
But this may be seen directly as
well: Let $\xi \in T_x(G(d,V))^u$ be an imaginary tangent of
weight $h\delta$, say. Weight considerations then yield $\xi e_{k}
= e_{k+hn}$, for $e_k \in x$, if $\xi e_k \neq 0$. Of course, $x$
contains $e_{k+hn}$ if $h \geq 0$, thus necessarily $h < 0$ if
$\xi \neq 0$.
\begin{defn}
Let $x \leq \ws \in I^u$. For any $i$ between $1$ and $n$ and $h >
0$ let $\xi_{i,h} \in T_x(G(d,V))$ be def\/ined as follows, provided
$h \leq \ell_i(x)$ and $h_i(x) - hn > 0$:
\[
\xi_{i,h}e_k = \begin{cases} p_x e_{k-hn}, \, & k \in S_i(x),\\ 0, \, &
\text{otherwise}.
\end{cases}
\] Similar to the real case, we simply
say $\xi_{i,h}$ is def\/ined at $x$, if $h_i(x) - hn > 0$.
\end{defn}

\begin{rem}
It is clear that $\xi_{i,h}$ is $u$-invariant. It is also clear,
that every imaginary tangent~$\xi$ of weight $-h\delta$ is a
linear combination of those $\xi_{i,h}$ which are def\/ined at~$x$:
As remarked above, $\xi$ is $u$-invariant and determined by its
values on the various $e_{h_i(x)}$, and $\xi e_{h_i(x)} =
\lambda_i e_{h_i(x) - hn}$, if $h_i(x) - hn > 0$, and zero
otherwise. Also, $\xi e_{h_i(x)} = 0$ if $h > \ell_i(x)$, for then
$\tau^{\ell_i(x)}\xi e_{h_i(x)} \neq 0$. Thus, $\xi = \sum_i
\lambda_i \xi_{i,h}$ provided we def\/ined $\lambda_i \xi_{i,h}$ to
be zero if $h_i(x) - hn \leq 0$ or $\ell_i(x) < h$.
\end{rem}
Our goal now is to describe those imaginary tangents which
actually appear in $T_x(X(\ws))$. To this end, let us keep $h > 0$
f\/ixed throughout the remainder of this subsection, and put
\begin{equation}\label{eq:SHX} S(h) = S(h,x)
= \{ i \mid h_i(x) - hn > 0; h \leq \ell_i(x) \},
\end{equation} the set of
indices for which $\xi_{i,h}$ is def\/ined. Clearly,
\begin{equation} \label{eq:GRASS_IMAGINARY}T_x(G(d,V))
^u_{-h\delta} = \bigoplus_{i \in S(h)} \K \xi_{i,h}
\end{equation}
is the weight space of imaginary weight $-h\delta$. Before we
describe the tangents belonging to $T_x(X(\ws))$ we need some more
notation. Recall the notion of Pl\"ucker coordinates on $G(d,V)$:
for $x \in I$, let $p_x$ be the corresponding Pl\"ucker coordinate
(equal to $e_{x_1}^* \wedge e_{x_2}^* \wedge \dots \wedge
e_{x_d}^*$, with the $e_{i}^*$ being a dual basis for the $e_i$).
The $p_x$ generate the homogenous coordinate ring of $G(d,V)$ (for
the Pl\"ucker embedding $G(d,V)\hookrightarrow\mathbb P
(\bigwedge^d V))$. For $x \in I$, let $U_x$ denote the open set
$p_x \neq 0$ in $G(d,V)$. Then $U_x$ is open af\/f\/ine, $\TV$-stable
(and actually isomorphic to a $\TV$-module). For any $\theta \in
I$, $f_\theta = \frac{p_\theta }{p_x}$ is a well def\/ined function
on $U_x$; $\mathcal O(U_x)$ is generated by those $f_\theta$ which
have nonzero dif\/ferential at $x$. The $\theta$'s in $I$ for which
this holds are precisely those, which dif\/fer in exactly one entry
from $x$, that is, $\theta = r_{i,j}x$ for suitable~$i$,~$j$. Here
$r_{i,j} \in \WV = \Sc_{dn}$ denotes the ref\/lection exchanging $i$
and $j$. For any $E_{kl} \in T_x(G(d,V))$ we have
$df_{\theta,x}(E_{kl}) = (-1)^{d(\theta)} \delta_{ik}\delta_{jl}$,
where $d(\theta)$ is the dif\/ference in positions between $j$ in
$x$ and $i$ in $\theta$. Let $\theta_{i,r}$ be obtained from $x$
by replacing $h_i(x) + rn$ with $h_i(x) + (r-h)n$, provided $0
\leq r < h\leq \ell_i(x)$ and $h_i(x)
 + (j- h)n > 0$.
\begin{lem}\label{lem:TRACE_RELATION}On $T_x(X(\ws))$ we have \[ \sum_{i\in S(h,x)}
(-1)^{d(\theta_{i,0})}df_{\theta_{i,0},x} = 0.\]
\end{lem}
\begin{proof}
Consider the action of $\tau^h$ on $\bigwedge^d V$ and
$\bigwedge^d V^*$ (acting as
\[v_1 \wedge v_2 \wedge \dots \wedge v_d \mapsto
\sum_{i=1}^d v_1 \wedge \dots \wedge \tau^h(v_i) \wedge \dots
\wedge v_d\] and similarly for $\bigwedge^d V^*$). If $M \subset
V$ is a $d$-dimensional subspace normalized by $u$, $\bigwedge^d
M$ is in the kernel of $\tau^h$ (as $\tau^h$ acts nilpotently on
$M$). Thus, for every $w \in I$, $p_w \in \bigwedge^d V^*$
satisf\/ies $p_w \tau^h(M) = 0$, and consequently $\tau^h p_w$
vanishes on $G(d,V)^u$. In particular $\tau^h p_x = 0$ on
$G(d,V)^u$. A straight forward computation (keeping in mind that
$\tau^h e^*_{k} = -e^*_{k-hn}$) shows that
 \begin{equation}\label{eq:TRACE_RELATION}
\tau^h p_x = \sum_{i = 1}^n\sum_{\substack{0 \leq j < \ell_i(x)\\
0 < h_i(x) + (j-h)n <h_i(x)}}
(-1)^{d(\theta_{i,j})}p_{\theta_{i,j}} + R,
\end{equation}
where $R$ is a linear combination of Pl\"ucker coordinates $p_w$
where $w$ dif\/fers from $x$ in strictly more than one element.
Localizing to $x$, we therefore obtain a relation
\[ \sum_i \sum_{\substack{0 \leq j < \ell_i(x) \\ 0 < h_i(x) + (j-h)n <
h_i(x)}} df_{\theta_{i,j},x} = 0\] on $T_x(X(\ws))$. Every summand
of this relation has $\hT$-weight $-h\delta$. To prove the lemma
it therefore suf\/f\/ices to consider the relation evaluated on
$T_x(G(d,V))^u_{-h\delta}$. Here, though, $df_{\theta_{i,j},x} =
(-1)^{d(\theta_{i,j})-d(\theta_{i,0})}df_{\theta_{i,0},x}$ by the
def\/inition of $\xi_{i,h}$, if $0 < h_i(x) - hn$. If on the other
hand $h_i(x) - hn = 0$, then $T_x(G(d,V))^u_{-h\delta}$ contains
no element supported in $E_{k,i+rn}$ for all $k$ and $r$, and
therefore $df_{\theta_{i,j},x} = 0$ in $T_x(G(d,V))^u_{-h\delta}$.
Summarizing, we get $\sum\limits_{i=1}^n
h(-1)^{d(\theta_{i,0})}df_{\theta_{i,0},x} = 0$, and the result
follows.
\end{proof}
Notice that $df_{\theta_{i,0},x}(\xi_{i,h}) =
(-1)^{d(\theta_{i,0})}$. Thus, $\sum\limits_{i\in S(h)}c_i \xi_{i,h} \in
T_x(X(\ws))$ implies $\sum_i c_i = 0$. We refer to the relations
of Lemma~\ref{lem:TRACE_RELATION} as \emph{trace relations}. The
reason for this is that when consi\-de\-ring the open immersion of the
nullcone $\mathcal N$ of nilpotent matrices into $X(w^1)$ alluded
to in the Introduction (cf.\ Lusztig's isomorphism \cite{LU2}),
$0$ is sent to $e$, and these relation in case $h = 1$ actually
correspond to the vanishing of the trace on nilpotent matrices.

We will show below  (cf.\ Theorem~\ref{thm:ONE_STRING_TANGENT}),
that the trace relations are the \emph{only} linear relations on
$T_x(X(w^s))_{-h\delta} \subset T_x(G(d,V))^u_{-h\delta}$.

Let $x,w\in I^u$, $x \leq w $. Similar to the case of real roots
(cf.\ Remark~\ref{rem:REALS_IN_Y}), we have the following Lemma
describing $T_x(Y(w))^u_{-h\delta}$. First one notation: Set
\begin{equation}\label{eq:SHXW}
S(h,x,w) = \{i \in S(h,x) \mid r_{h_i(x)+(j-h)n,h_i(x) + jn}x \leq
w,\, \forall \, j = 0,1,\dots,h-1\}.
\end{equation} (Here, $S(h,x)$ is as in equation~\eqref{eq:SHX}.)

\begin{lem}\label{lem:IMAGINARY_TANGENTS} Let $x,w\in I^u$, $x \leq w $. Then
\begin{enumerate}\itemsep=0pt
\item[$(1)$] $T_x(Y(w))^u_{-h\delta}$ is spanned by all $\xi_{i,h}$,
defined at $x$, which are contained in $T_x(Y(w))$.
\item[$(2)$]
$T_x(Y(w))_{-h\delta}^u = \bigoplus_{i \in S(h,x,w)} \K
\xi_{i,h}$.
\end{enumerate}
\end{lem}
\begin{proof}
Indeed, any $\xi \in T_x(Y(w))^u_{-h\delta}$ is a linear
combination of $\xi_{i,h}$s, def\/ined at $x$. On the other hand,
the $\xi_{i,h}$s are supported in entirely dif\/ferent
$\TV$-eigenspaces. Thus, $\xi\in T_x(Y(w))^u$ can be supported in $\xi_{i,h}$ only if all the $\TV$-eigenvectors in which $\xi_{i,h}$ has a nonzero component -- and thus $\xi_{i,h}$ -- are elements of $T_x(Y)$.
Assertion (1) now follows.

Assertion (2) follows from Assertion (1) and the def\/inition of
$S(h,x,w)$.
\end{proof}

\begin{figure}[t]
\centerline{\includegraphics{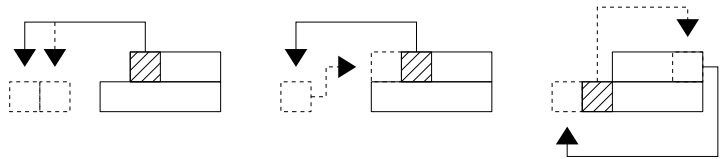}}
\caption{An example for how to create
imaginary tangents as in the proof of
Lemma~\ref{lem:IMAGINARY_SING}. The left picture shows part of
$\Sigma(w)$ and the ef\/fect of $s_\halpha$ and $s_\hbeta$ on the
f\/irst box of the row (the arrows with the solid lines refer to
$s_\halpha$). The picture in the middle shows the ef\/fect of
$\xi_{-\beta}\xi_{\halpha}$ and the picture to the right shows
$\xi_{\halpha}\xi_{-\hbeta}$.}\label{fig:IMAGINARY}
\end{figure}

As a consequence, we have the following
\begin{thm}\label{kap1} $T_x(Y(\ws))^u
= T_x(G(d,V))^u$ for all $x \in I^u$.
\end{thm}
 \begin{proof} It is easily seen that
for all $x$, $S(h,x,\ws) = S(h,x)$.  The result follows from this
fact, the above lemma, equation~\eqref{eq:GRASS_IMAGINARY}, and
Corollary~\ref{cor:CURVE_COR}.
\end{proof}
Returning to our study of $T_x(X(\ws))$ for some arbitrary but
f\/ixed $x \leq \ws$, we make the following
\begin{defn}Let $i,j \in S(h,x)$ be arbitrary but distinct. For
$\alpha = (ij) \in \Phi$ we set
\[ \xi_{\alpha,h} = \xi_{i,h} - \xi_{j,h} \in T_x(G(d,V))^u.\]
\end{defn}
Notice that by the trace relation, the tangents of the form
$\xi_{\alpha,h}$ span $T_x(X(w))_{-h\delta}$.
\begin{lem}\label{lem:IMAGINARY_SING}
Suppose $s_\halpha$ is a large reflection defined at $x \leq w \in
I^u$ with $x < s_\halpha x \leq w$. Then~$x$ is a singular point
of $X(w)$.
\end{lem}
\begin{proof}
We will show that there is an imaginary tangent in $T_x(X(w))$. In
fact, it will be $\xi_{\Re(\halpha),1}$ (recall $\Re(\halpha)$
from Section~\ref{weyl}). By Lemma~\ref{lem:IMAGINARY_TANGENT}, the
conclusion follows.

We may write $\halpha = (ij) + h\delta$, with $h \leq -1$. Let
$\hbeta = \halpha + \delta$. $\hbeta$ is negative as well: if
$(ij) < 0$ this is clear; if $(ij) > 0$ then $h \leq -2$ since
$s_\halpha$ is large, and $\hbeta < 0$ follows. There are two
cases to consider: either $s_\hbeta x$ is def\/ined at $x$, and
$s_\hbeta x > x$, or $s_\hbeta x = x$. In both cases $U_{-\hbeta}$
is a subgroup of the stabilizer of $x$ in $\calB$ and consequently
acts on $T_x(X(w))$. Notice that the Lie algebra of~$U_{-\hbeta}$ is spanned by
$\xi_{-\hbeta} \in \End(V)$ (cf.\ equation~\eqref{eq:END_XI}). A
straight forward computation then shows that
$[\xi_{-\hbeta},\xi_{\halpha}] = -\xi_{(ij),1}$ (see
Fig.~\ref{fig:IMAGINARY}). As this is the action of $\Lie
U_{-\hbeta}$ on $T_x(X(w))$, it follows that $\xi_{(ij),1}\in
T_x(X(w))$.

Of course the lemma also follows immediately if one considers the
fact that $[\hat \fg_{-\hbeta},\hat \fg_{\halpha}] \subset \hat
\fg_{-\delta}$ is nonzero and therefore a tangent of $T_{x}(X(w))$
(since $\hat \fg \to T_x(\GL_n(F)/\GL_n(A))$ is surjective with
kernel $x \hat \fg_+ x^{-1}$ where $\hat \fg_{+} =
\bigoplus_{\delta(\halpha) \geq 0} \hat \fg_\halpha$).
\end{proof}

\begin{rem}\label{rem:Peterson}
Recall (cf.~\cite{LSe}) that in the classical setting in type $A$,
a point $x \leq w$ is smooth in $X(w)$ if and only if there are
precisely $\dim X(w)$ ref\/lections $r$ such that $x\neq rx \leq w$. In
the af\/f\/ine setting, this is no longer true; this description fails
for example if one of these ref\/lections is large.

However, using  D. Peterson's ideas of deforming tangent spaces
(see \cite{CK} for a discussion of this approach), it seems to be possible to show that if for all
$y$ with $x \leq y \leq w$ we have $TE(X(w),y)$ has $\dim X(w)$
elements, and furthermore no ref\/lection $r$ with $y < ry \leq w$
is large, then $x$ is a~smooth point.
\end{rem}

\section{Real and imaginary patterns}\label{sec:PATTERNS}
We are now ready to describe several types of singularities of a
given $X(w) \subset X(\ws)$. As it turns out, the singularities
are best described using $L(w)$, due to the subtlety that when
down-exchanging $i > j$ (with $\ell_i(w)> \ell_j(w)$) the result
is not just ex-changing the rows in $\Sigma(w)$.

\subsection{Imaginary patterns} The previous section of course provides a very elementary
way of producing singularities. For the sake of consistency we
give it a name:
\begin{defn}
An \emph{imaginary pattern} $P$ in $L(w)$ is a pair of integers
$(i,j)$ ($i < j$) with $i \leq n$, such that $L(w)_i > L(w)_j +
1$. For such a pattern $P$, let $w_P$ be obtained from $w$ by
repla\-cing~$h_i(w)$ with $h_j(w)-n$.
\end{defn}
\begin{rem}
Notice that $w_P$ is clearly singular, because it is of the form
$sw < w$ with $s$ a large ref\/lection.
\end{rem}
In some cases all maximal singularities of a given Schubert
variety arise in this fashion; for instance the single maximal
singularity of $X(\ws)$ is $\ws_P$ for $P = (12)$ (see
Section~\ref{sec:ONESTRING}).

Obviously, the condition of not admitting \emph{any} imaginary
pattern forms a serious obstruction against being
non-singular. It is immediately forced that for a smooth
$X(w)$, $\ell_i(w) \leq \ell_j(w) + 2$ for all pairs $i$, $j$, and
$\ell_i(w) = \ell_j(w) + 2$ is possible only if $j < i$.

\subsection{Real patterns}

Perhaps more interesting are the singularities which arise
because $TE(X(w),x)$ is too large, or, in other words, because
there are too many $\hT$-stable curves through $x$. Recall from
Section~\ref{sec:REFLECTIONS} that the singularities of $Y(w)$
correspond to the hooks in $\CM(w)$. This is no longer true for
af\/f\/ine Schubert varieties, but there is a type of pattern in
$L(w)$ which closely resembles this concept. In fact, a hook in
$\CM(w)$ is more or less a ``gap'' in $w$, i.e.\ an index $i$ such
that $w_{i+1} \neq w_i + 1$ together with a (f\/irst) position $k >
i$ such that $w_{i+1} + (k-i) \neq w_{k+1}$. For $X(w)$ this is
more complicated:
\begin{defn}
Let $w \in I^u$. A \emph{real pattern of the first kind} in $L(w)$
is a sequence of integers $1\leq i < g < j < k \leq 2n$ subject to
the following conditions:
\begin{enumerate}\itemsep=0pt
    \item[$1)$] $i < n$, $j < i + n$, and $k < g + n$;
    \item[$2)$] $l_i \geq l_j > l_g \geq l_k$.
\end{enumerate}
If $P = (i,g,j,k)$ is such a pattern, $w_P$ is obtained from $w$
by the following sequence of down-exchanges: $i \da g$, $g \da j$,
$g \da k$.
\end{defn}

\begin{figure}[t]
\centerline{\includegraphics{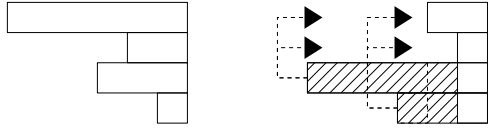}}
\caption{The basic example of a real
pattern of the f\/irst kind. Here $P = (1,2,3,4)$, and $w_P$ is
shown on the right, together with four up-exchanges. Notice that
according to Lemma~\ref{lem:CODIM}, $w_P$ has codimension~$3$.}\label{fig:PATTERN_1ST}
\end{figure}

\begin{prop}\label{prop:PATTERN_FIRST} If $P = (i,g,j,k)$ is a real pattern of the first kind
in $L(w)$, then $w_P$ is singular. More precisely, if $k < i + n$,
then $\vert E(X(w),w_P)\vert > \dim X(w)$. If $k > i + n$, then
$w_P$ admits a large reflection $s$ such that $w_P < sw_P \leq w$.
\end{prop}
It is worth mentioning, that $k = i + n$ does not occur, because
$L(w)_i > L(w)_k$, which never holds for $k = i + n$. Before
proving Proposition~\ref{prop:PATTERN_FIRST}, notice f\/irst that
\begin{rem}\label{rem:LW_P}
$w_P$ is determined by the requirements $L(w_P)_i = L(w)_g$,
$L(w_P)_g = L(w)_t$, $L(w_P)_j = L(w)_i$, and f\/inally $L(w_P)_t =
L(w)_j$.

This is an immediate consequence of the fact, that for each
down-exchange in the def\/inition of $w_P$ the corresponding indices
$i',j'$ satisfy $i' <j'$ and $L(w')_{i'} > L(w')_{j'}$ where $w'$
denotes the intermediate step on which the down-exchange is
performed.
\end{rem}
It  should be mentioned that it is possible that Gasharov's proof of similar statements in the classical case \cite{gash} could be adapted to our situation to simplify the proofs of Proposition~\ref{prop:PATTERN_FIRST} as well as Proposition~\ref{prop:PATTERNSND}. However, except as mentioned below in Remark~\ref{rem:PATTERN}, we don't see how.

\begin{proof}[Proof of Proposition~\ref{prop:PATTERN_FIRST}]
Let $L(w) = (l_1,l_2,\dots,l_{2n})$. First suppose that $k < i +
n$. We will show that the number of ref\/lections $s$ with $w_P <
sw_P \leq w$ is too big, i.e.\ strictly larger than the
codimension of $X(w_P)$ in $X(w)$ (see
Lemma~\ref{lem:TOOMANYCURVES}).

Let $w_1$ be obtained from $w$ by down-exchanging $i$ and $g$
(cf.\ Fig.~\ref{fig:REAL_FIRST_REFLECTIONS}). Then $w > w_1 >
w_P$. Let~$c$ be the codimension of $w_P$ in $X(w_1)$. By
Deodhar's Inequality there are at least~$c$ ref\/lections~$s$ such
that $w_P < sw_P \leq w_1$. Let~$c_1$ be the codimension of
$X(w_1)$ in $X(w)$. We will construct $c_1 + 1$ ref\/lections $s$
with $w_P < sw_P \leq w$ but $sw_P \nleq w_1$.

\begin{figure}[t]
\centerline{\includegraphics{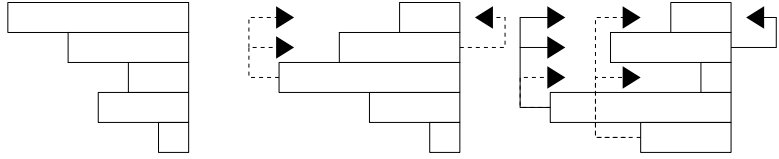}}
\caption{Another real pattern of
the f\/irst kind: left $w$, in the middle $w_1$, and to the right
$w_P$; here $P = (1,3,4,5)$. The up-exchanges constructed from the
ones corresponding to $w_1 < w$ have a solid line.}\label{fig:REAL_FIRST_REFLECTIONS}
\end{figure}

Consider the $c_1$ ref\/lections corresponding to the down-exchange
of $i$ and $g$ in $w$; such a~ref\/lection $s$ satisf\/ies  $w \geq
sw_1> w_1$. According to
Remark~\ref{rem:CORRESPONDING_REFLECTIONS} there are three kinds
of these ref\/lections: f\/irst, an up-exchange of an element $h \in
g_\geq(i,g,w)$ with $i$, and second an up-exchange of such an
element with $[g]$. Finally, there is also the up-exchange of~$i$
and $[g]$ (turning $w_1$ into~$w$).

Now consider the f\/irst type. i.e.\ $s$ is an up-exchange of some
$h \in g_\geq(i,g,w)$ with $i$, and then $w_1 < sw_1 < w$. Also,
$i < h < g$ and $L(w)_i \geq L(w)_h > L(w)_g$. Notice that this
means $[h] \neq [k]$; $[h] \neq [j]$ as well, because $i < h < j <
i + n$. It follows that the very same up-exchange is def\/ined at
$w_P$ and $sw_P > w_P$. On the other hand $sw_P \leq w$, because
it may be obtained by $g \da j$, and then $g \da k$ applied to
$sw_1$. Also $sw_P \nleq w_1$: of course, $w$, $w_1$, $w_P$, $sw_P$ all
dif\/fer only in the $n$-strings through $i$, $h$, $g$, $j$, $k$ (in fact, this
means, we may actually assume $n = 5$); elements of these strings
will be referred to as \emph{relevant}. Moreover, $L(sw_P)_i =
l_h$, so the elements of $sw_P$ in these strings that are less
than or equal to $h_i(sw_P)$ are given by $h_i(sw_P)$ itself, the
f\/irst $l_i - l_h$ elements of $S_{[j]}(sw_P)$, and the f\/irst $l_j
- l_h$ elements of $S_{[k]}(sw_P)$ (if $l_j > l_h$). Notice that
if for any integer $r$ among $h$, $g$, $j$, $k$ we have $r > n$, then $[r]
< i$ by assumption. So to see why $l_i - l_h$ elements of
$S_{[j]}(sw_P)$ are less than or equal to $h_i(sw_P)$, observe
that if $j \leq n$, then this is clear, as then also $h \leq n$.
Otherwise, it follows that $[j] < i$ as $j < i + n$. Thus, exactly
$l_h - 1$ entries of $S_{[j]}(sw_P)$ are strictly larger than
$h_i(sw_P)$. On the other hand, $L(sw_P)_{[j]} = l_i - 1$, so the
dif\/ference, i.e.\ the number of those less than or equal to
$h_i(sw_P)$ is precisely $l_i - l_h$. The case of~$[k]$ is
similar. Notice that $k < i + n$ is needed here only if $l_j >
l_h$.

In $w_1$, however, the relevant elements less than or equal to
$h_i(sw_P)$ comprise only $l_i-l_h$ elements of $S_{[g]}(w_1)$ and
possibly $l_j - l_h$ elements of $S_{[j]}(w_1)$ (if $l_j > l_h$).
Thus the total number is strictly smaller, and $sw_P \nleq w_1$.

The second possibility is that $s$ is an up-exchange of $h$ and
$g$. In this case, let $s'$ be the ref\/lection associated to $j \up
h$  in $w_P$. As $L(w_P)_h < l_i$ in this case, $s'$ is
well-def\/ined. Notice that this case also includes $h = i$.
Moreover $s'w_P \leq w$, as it may be obtained by $j \da k$ and
then $g \da j$ in $sw_1$. Again $sw_P \nleq w_1$, since the number
of elements less than or equal to $h_{[h]}(sw_P)$ is strictly
larger than the same number for $w_1$.

Summarizing, each of the $c_1$ ref\/lections at $w_1$ gives rise to
a ref\/lection $s$ at $w_P$ such that $w_P < sw_P \leq w$, but $sw_P
\nleq w_1$.

It remains to construct one additional ref\/lection with this
property. Let $s_1$ be the ref\/lection corresponding to $k \up i$ in
$w_P$. Clearly $s_1 w_P \leq w$ as $s_1w_P$ may be obtained from $w$ by
down-exchanging $i$ and $j$, and $g$ and $k$ (which uses $k < i +
n$). But again, $s_1w_P \nleq w_1$: the number of elements in
$s_1w_P$ less than or equal to $h_i(s_1w_P)$ is strictly larger
than the same number computed for $w_1$. In $w_1$ the only
relevant such elements are the f\/irst $l_i - l_j$ elements of
$S_{[g]}(w_1)$. The same number for $s_1 w_P$, however, is given
by $1$ for $h_i(s_1w_P)$, plus the f\/irst $l_i -l_j$ elements of
$S_{[j]}(s_1w_P)$.

Finally, in the case $k > i +  n$, we will construct a large
ref\/lection: let $s$ correspond to two subsequent up-exchanges of
$k$ and $i$ in $w_P$. Notice that $s$ is indeed large. $sw_P$ is
obtained from~$w_P$ by decreasing $L(w_P)_k$ by one, and
increasing $L(w_P)_i$ by one ($k > n$, hence $l_k > 0$ and
therefore $L(w_P)_i = l_g \geq l_k > 0$). An example is outlined
in Fig.~\ref{fig:REAL_FIRST_LARGE}.

\begin{figure}[t]
\centerline{\includegraphics{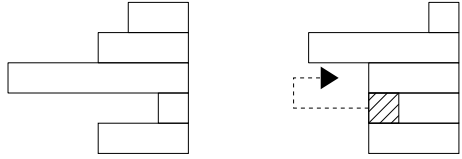}}
\caption{A real pattern of the f\/irst
kind with $k > i + n$. $w$ is shown to the left, and $w_P$,
together with the large up-exchange, is shown to the right.
Assuming $n = 5$, $P = (3,6,7,9)$.}\label{fig:REAL_FIRST_LARGE}
\end{figure}

To recap, $sw_P$ is characterized by $L(sw_P)_i = l_g + 1$,
$L(sw_P)_g = L(w_P)_g = l_k$, $L(sw_P)_j = L(w_P)_j = l_i$, and
$L(sw_P)_k = l_j - 1$. Thus, $sw_P \leq w$, as it is obtained from
$w$ by $i \da j$, and $g \da k$, and, if $l_j > l_g + 1$, $i \da
k$. As $s$ is large, $w_P$ must be a singular point of $X(w)$, as
claimed.
\end{proof}

The real patterns of the f\/irst kind are modeled loosely after the
hooks in the classical setting, $h_{[g]}(w)-n$ playing the role of
the ``gap'' between elements of $S_i(w)$ and $S_{[j]}(w)$.
There is
another kind of pattern for which this analogy fails:
\begin{defn}
Let $w \in I^u$. A \emph{real pattern of the second kind} for $w$
is a sequence of integers $i < j < g < k$ subject to the following
conditions:{\samepage
\begin{enumerate}\itemsep=0pt
\item[$1)$] $i \leq n$; $g < i + n$; $k < j + n$;
\item[$2)$] $L(w)_j > L(w)_i
\geq  L(w)_k > L(w)_g$.
\end{enumerate}
If $P = (i,j,g,k)$ is such a pattern, $w_P$ is obtained from $w$
by $j \da k$, $i \da j$, and f\/inally, $i \da g$.}
\end{defn}

Again, $k$ is never equal to $i + n$ in such a pattern. However,
the case $k = i + n$ will be what we call an exceptional pattern
below. $w_P$ is def\/ined by $L(w_P)_i = L(w_P)_g$, $L(w_P)_j  =
L(w_P)_i$, $L(w_P)_g = L(w_P)_k$, and $L(w_P)_k = L(w_P)_j$.

\begin{figure}[t]
\centerline{\includegraphics{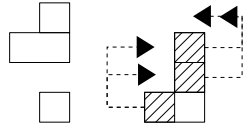}}
\caption{A real pattern of the second kind
(left). Here $n = d =4$, and the pattern is $(1,2,3,4)$. $w_P$ is
shown in the middle, together with the four up-exchanges def\/ined
at $w_P$ (note that $\codim_w(w_P) = 3$).}\label{fig:REAL_SECOND}
\end{figure}

\begin{prop}\label{prop:PATTERNSND} If $P$ is a real pattern of the second kind for $w \in
I^u$, then $w_P$ is singular in~$X(w)$. In fact, $\vert
E(X(w),w_P) \vert > \dim X(w)$, if $i + n > k$.
\end{prop}
\begin{proof}
The reasoning is similar to the case of real patterns of the f\/irst
kind. Again, we f\/irst assume that $i + n > k$. Let $w_1$ be
obtained from $w$ by down-exchanging $j$ and $k$. Then $w > w_1 >
w_P$. Let $c_2$ be the codimension of $X(w_P)$ in $X(w_1)$. Then
there are at least $c_2$ up-exchanges $s$ at $w_P$ such that $sw_P
\leq w_1$.

Let $c_1$ be the codimension of $X(w_1) \subset X(w)$; then there
are $c_1$ corresponding ref\/lections, all of them involving $j$ or
$k$ (and exactly one, both). For each such ref\/lection $s$, let
$s'$ be the up-exchange of $w_P$ as follows: if $s$ is $k \up j$,
then $s = s'$ is def\/ined. Suppose $s$ involves $h \in
G_{\geq}(w,j,k)$: if $s$ up-exchanges $k$ and $h$, then the very
same up-exchange is def\/ined at $w_P$, and $s' = s$ (as then $l_j =
L(w_P)_k > L(w_P)_h$). Otherwise, $s$ exchanges $h$ and $j$ (and
then $l_h > l_k$); if $l_h > l_i$, the same up-exchange is def\/ined
and again $s = s'$. Finally, if $l_h \leq l_i$, then if $h > g$,
replace $s$ by $h \up g$; otherwise, if $h < g$, by $h \up i$.
Notice that the case $h = g$ does not occur, because $l_g < l_k$.

In all these cases, $s'w_P \leq w$: If $s'$ is $j \up k$, then
$s'w_P$ is obtained from $w$ by $i \da k$ and then~$i\da g$. If
$s'$ involves $k$, but not $j$, $s'w_P$ obtained from $w$ by $j
\da h$, $j \da k$, $i \da j$, and $i\da g$. If~$s'$ is $h \up j$
for some $h$, then $s'w_P$ is the result of $j \da h$, $h \da k$,
$i \da g$, and $i \da h$. If $s'$ is $h \up i$, then recall that
$h < g$, and $l_h \leq l_i$; so $s'w_P$ is obtained from $w$ by $j
\da h$, $i \da j$, $j \da g$, and $g \da k$. If $s'$ is $h \up g$,
then $h > g$, and $l_h \leq l_i$; so $s'w_P$ is obtained from $w$
by $j \da g$, $i \da j$, $g \da h$, and f\/inally, $h \da k$.

Notice that none of the $s'$ are among the $c_2$ up-exchanges
corresponding to $w_1 > w_P$ since $s'w_P \nleq w_1$: this is
clear if $s'$ involves $k$, as $l_j = L(w_P)_k$ is the largest
relevant length. So suppose otherwise; if $s'$ involves $h$ and
$j$, then $l_h > l_i$, and so the position of $h_j(s'w_P)$ in
$sw'_P$ is equal to the position of $h_h(w_1)> h_j(sw'_P)$ in
$w_1$. If $s'$ involves $h$ and $i$, then an easy calculation
shows that
\[ \vert \{ m \in w_1 \mid m \leq h_i(s'w_P) \}\vert < \vert \{ m \in s'w_P \mid m \leq h_i(s'w_P)\} \vert. \]
Similarly, if $s'$ is the up-exchange of $h$ and $g$, then the
same is true with $h_i(s'w_P)$ replaced by~$h_g(s'w_P)$. Thus, we
have a combined total of $c_1 + c_2$ up-exchanges. But there is at
least one additional up-exchange, namely the one of $k$ and $g$:
Notice that since $g \not \in G_{\geq}(w,j,k)$, this one is not
among the $c_1$ ref\/lections $s'$. Moreover, the result is not
contained in $w_1$ since $L(w_P)_k = l_j$ is the longest relevant
length and would have to be at the same position as in $w_1$. All
in all, this shows that $\dim E(X(w),w_P) > \dim X(w)$, and we are
done.

Now suppose that $k > i + n$. Then $i < [k]$ and $l_i > l_{[k]}$,
and it is possible, that the up-exchange of $g$ and $k$ in the
reasoning above has been listed before. However, we may
up-exchange $g$ and~$j$ twice, resulting in $w' \leq w$,
satisfying $L(w')_i = l_g$, $L(w')_j = l_i + 1$, $L(w')_{k} =
l_j$, $L(w')_g = l_{k} - 1$ (this last value is $\geq 1$ if $g >
n$, as then $l_k > l_g \geq 1$). $w'$ is obtained from $w$ by
down-exchanging~$j$ and~$i$ (if $l_j > l_i + 1$), and then
down-exchanging~$i$ and~$k$, and then $i$ and $g$. As a double
up-exchange it corresponds to a large ref\/lection, $w_P$ is
singular.
\end{proof}

\begin{rem}\label{rem:PATTERN}
As pointed out by one of the referees, in some instances the fact that $w_P$ singular can be seen quicklier. The situation is as follows: Let $L(w) = (l_1,l_2,\dots,l_{2n})$ and suppose $P = (i,g,j,k)$ is a real pattern of the f\/irst kind, say, where $k < i + n$. For $\sigma \in \Sc_n$, the symmetric group in $n$ letters, viewed as the permutations of $\{ i.i + 1,\dots,i+n\}$, def\/ine $\sigma w$ as the element obtained from $w$, by permuting $(i,i+1,\dots,i + n-1)$ according to $\sigma$.
Let $L_1 \leq L_2 \leq \dots \leq L_k$ be the \emph{distinct} values of $l_i,l_{i+1},\dots,l_{i+n-1}$, and f\/inally put $d_j = \vert \{ i \leq t <  i + n \mid \ell_t(w) = L_i \}\vert$. Consider the variety $\calF = \calF(d_1,d_2,\dots,d_k)$ of partial f\/lags $0 \subset V_1 \subset V_2 \subset \dots \subset V_k \subset \CC^n$ with $\dim V_i = d_i$. $w$ naturally def\/ines a point in $\calF$, which in turn can be described by any element $\tau$ in $\Sc_n$, for which $\tau(i) ,\tau(i+1), \dots , \tau(i + d_1 -1)$ are the indices of the basis elements in $V_{1}$ (here, the indices of the $l_i$ for which $l_i = L_1$), and so on. If we choose $\tau$ correctly, the existence of the pattern $P$ then means that $\tau$ may be chosen in a way such that  $i$, $g$, $j$, $k$ appears in order $k$, $g$, $j$, $i$ which is a (Type~II) pattern in the classical sense. In particular, $w_P$ then corresponds to $\tau_P$ obtained from $\tau$ by reordering $k$, $g$, $j$, $i$ to $g$, $i$, $k$, $j$ (and indeed, $w_P$ is obtained by replacing $l_i$ with $l_g$, $l_j$ with $l_i$, $l_g$ with $l_k$, and $l_k$ with $l_j$.

Also, if $r$ is a transposition exchanging $p$ and $q$, say, then $\tau_P < r\tau_P \leq w_P$ if and only if $w_P < rw_P \leq w$ (if we choose the right order on $\Sc_n$). Here $rw_p$ means switching $l_p$ and $l_q$ (corresponding to an up- or down-exchange, or if $l_p = l_q$ to doing nothing). Since $\tau_P \leq \tau$ is a~singular point of the corresponding classical Schubert variety, $w_P$ is singular in~$w$. It is very likely that this can carried over to some of the other patterns. However, the situation is unclear when $i + n < k$ and it seems that we cannot avoid these patterns to f\/ind ``combinatorial'' singularities (i.e.\ points where we have too many $\hT$-stable curves); see Remark~\ref{rem:SMOOTH} for an example where an exceptional pattern (def\/ined below) is needed.
\end{rem}

\subsection{Exceptional patterns} We conclude this section with two
degenerations of real patterns.

\begin{defn}
Let $w \in I ^u$ and $L(w) = (l_1,\dots,l_{2n})$. Then a sequence
$P= (i,g,j=i+n,k)$ with $i < g < j < k$ is called an
\emph{exceptional pattern of the first kind}, if $ k < g + n$, and
if $l_i > l_{[g]} + 1$, and $l_g \geq l_k$. The associated point
$w_P$ is def\/ined by f\/irst down-exchanging $i$ and $g$, and then
down-exchanging $g$ and $j$ (only if $l_i > l_g + 1$), and f\/inally
down-exchanging $g$ and $k$ as before.
\end{defn}
Notice that if $g > n$, and $l_i = l_{[g]} + 1$, the second
down-exchange of $i$ and $g$ is void. Notice that the requirement
$k < g + n$ is not really necessary. If $k > g + n$, then $g < n$.
The above procedure then results in a point $w_P < w_Q$ for the
imaginary pattern $Q= (i,g)$, which clearly is singular. In fact
this immediately shows that $w_P$ is singular, whenever $g \leq
n$. In general, the argument is similar:
\begin{lem}\label{lem:SING_EXC_FIRST} Let $P$ be an exceptional pattern of
the first kind for $w$.  Then $w_P$ is a singular point of $X(w)$.
\end{lem}
\begin{proof}
By the remarks preceding the lemma, we may assume that $g > n$.
Then $w_P$ is characteri\-zed by $L(w_P)_{[g]} = l_{[k]}$, $L(w_P)_i =
l_i - 1$, $L(w_P)_{[k]} = L(w)_{[g]} + 1$. As $[k] > i$, and $l_i
\geq l_g + 1 = l_{[g]} + 2$, it follows that $w_P = w'_Q$ where
$w'$ is obtained from $w$ by $g \da k$, and $Q$ is the imaginary
pattern $(i,[k])$ for $w'$. Hence, the claim.
\end{proof}

\begin{figure}[t]
\centerline{\includegraphics{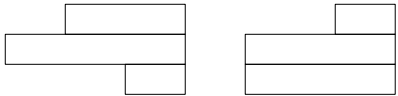}}
 \caption{$w$ (left) and $w_P$ for the exceptional pattern of the f\/irst kind $P = (2,4,5,6)$ (assuming $n = 3$).}\label{fig:EXCEPTIONAL1}
\end{figure}

\begin{rem}
 It is clear that an exceptional pattern of the f\/irst kind is only interesting if $g > n$, and $l_i = l_{g} + 1$. In all other cases, $w_P < w_Q$ where $Q = (i,g)$.
\end{rem}

As for the second kind, the situation is similar:
\begin{defn}
Let $w\in I ^u$, then a sequence $P = (i,j,g,k = i + n)$ is called
an \emph{exceptional pattern of the second kind} for $w$ or $L(w)
= (l_1,\dots,l_{2n})$, if $i < j < g < k$, and if $l_i < l_j$, and
$l_i > l_{g}$.

Then $w_{P,1}$ is def\/ined as $w'_Q$ where $w'$ is obtained by $i
\da j$, and $Q$ is the imaginary pattern $([j],[g])$ (resp.\
$([j],g)$) in $w'$, if $[j] < [g]$ (resp.\ $[j] > [g]$); $w_{P,2}$
is def\/ined as $i \da j$, $i \da j$, $i \da g$ applied to $w$.
\end{defn}
Notice that $Q$ is indeed an imaginary pattern in $w'$, because
$L(w')_i = l_j - 1$, $L(w')_j = l_i + 1$, $L(w')_g = l_g < l_i$.
Since $w'_Q$ admits a large ref\/lection relative to $w'$, $w'_Q$ is
singular in $X(w)$.

As for $w_{P,2}$, it need not always be singular. Indeed, if $l_j
= l_i + 1$, then $w_{P,2}$ is just $i \da g$, $j \da g$ applied to
$w$. However, if $l_i < l_j - 1$, then $w_{P,2} = w''_Q$ where
$w''$ is obtained by $i \da g$, and $Q = (j,g)$. Also notice that
$w_{P,2}$ is the point obtained by applying the rule for secondary
kind patterns (ignoring that $k = i + n$).

\begin{rem}
The exceptional patterns of the second kind are interesting only,
if $l_g = l_i - 1$. For suppose $l_g < l_i - 1$. Then, we actually
have $w_{P,1/2} \leq w_Q$ for $Q = (i,g)$. Also, if $l_i < l_j -
2$, then $w_{P,2} \leq w_Q$ for $Q = (i,j)$.
\end{rem}

\begin{figure}[t]
\centerline{\includegraphics{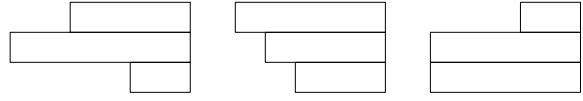}}
 \caption{$w$ (left), $w_{P,1}$, and $w_{P,2}$ (right) for the exceptional pattern $P = (1,2,3,4)$ of the second kind (assuming again that $n = 3$). Note that $w_{P,2}$ equals $w_Q$ where $Q = (2,4,5,6)$ is the pattern of the f\/irst kind from Fig.~\ref{fig:EXCEPTIONAL1}; it is often the case, that exceptional patterns of the f\/irst kind give rise to a pattern of the second kind.}\label{fig:EXCEPTIONAL2}
\end{figure}

\section{Two classes of Schubert varieties}\label{sec:ONESTRING}

We end this note by studying two classes of af\/f\/ine Schubert varieties where our discussion determines the singular locus completely.

\subsection[$\calP$-stable Schubert varieties]{$\boldsymbol{\calP}$-stable Schubert varieties}\label{subsec:PSTABLE}

Recall that an af\/f\/ine Schubert variety $X(w)$ is $\calP$-stable if
and only if $\ell_1(w) \geq \ell_2(w) \geq \dots \geq \ell_n(w)$: this is clearly necessary, for in order to be $\calP$-stable, it must be $\SL_n(\K)$-stable, and therefore invariant under the f\/inite Weyl group $W = W_\calP$; conversely being stable under $W_\calP$ is clearly enough by the Bruhat decomposition, and it is not hard to see that $X(w)$ is $W$-stable if and only if $Ww \subset X(w)$.

In other
words, $X(w)$ is $\calP$-stable if and only if $\Sigma(w)$ is an
actual Young diagram. It is worth mentioning, that $\calP$ acts
linearly on $V$. In fact, if $\calP_0 = \GL(A)$, then the image of
$\calP_0$ in $\GL(V)$ is precisely $C(u)$, the centralizer of $u$.
And, as far as the action on $G(d,V)$ is concerned, $\calP_0$~and~$\calP$ have the same orbits. By abuse of language we say $w \in
I^u$ is \emph{$\calP$-stable}, if $X(w)$ is. Since $e$ is the only
$\calP$-f\/ixed point in $X(\ws)$, this should not lead to confusion.

For $\calP$-stable elements of $I^u$, the Bruhat--Chevalley ordering
is much simpler: if $x,w \in I^u$ are $\calP$-stable, then $x \leq
w$ if and only if for $k = 1,2,\dots,n$:
\begin{equation}\label{eq:BRUHAT_PSTABLE}
\sum_{i = 1}^k \ell_i(x) \leq \sum_{i=1}^k \ell_i(w).
\end{equation}
Notice, that if $\Sigma(w)$ is a Young diagram, there are no real
patterns. This is clear for those of the f\/irst kind. For those of
the second, assume $P = (i,j,g,k)$ is such a pattern, then $L(w)_i <
L_(w)_j$ is possible only if $j > n$. But the requirement that
$L(w)_i \geq L(w)_k > L(w)_g$ then implies that $g > k$, a
contradiction. But there may be imaginary patterns: let
$1 \leq i_1 < i_2 < \dots < i_k < n$ be the uniquely determined
sequence of those integers satisfying $\ell_{i_k +1 }(w) <
\ell_{i_k}(w)$. If $w \neq e$, there is of course at least one such
integer. For each $r = 1,2,\dots,k$, let $j_r$ be the minimal index,
such that $\ell_{j_r}(w) < \ell_{i_r}(w) - 1$. Such an
index does not necessarily exist for the last integer $i_k$. If it doesn't, we remove
$i_k$ from the list and replace $k$ with $k-1$. It always exists for $i_{k-1}$. Thus, if $w \neq e$, there
is at least one pair $(i_r,j_r)$.
\begin{prop}\label{prop:PSTABLE}
Let $w \in I^u$ be $\calP$-stable. Keeping the notation above, for
$1 \leq r \leq k$, $P_r = (i_r,j_r)$ is an imaginary pattern of
$L(w)$, and the maximal singularities of $X(w)$ are among the points
$w_{P_r}$ \tu{(}$1 \leq r \leq k$\tu{)}. More precisely, $w_{P_r}$ is maximal
if and only if $i_r = \max \{ i_h \mid j_h = j_r\}$.
\end{prop}
\begin{proof}
It is clearly safe to assume that $w \neq e$. We have noted that
$\Sigma(w)$ is a Young diagram. As $X(w)$ is $\calP$-stable so is
its singular locus. Consequently, if $x$ is a maximal singularity of
$X(w)$, then $\Sigma(x)$ is a Young diagram as well. By
\eqref{eq:BRUHAT_PSTABLE} we have for all $h = 1,2,\dots,n$
\begin{equation}\label{eq:CRITERION} \sum_{j=1}^h \ell_j(x) \leq
\sum_{j=1}^h \ell_j(x).
\end{equation}
 Let $i$ be the f\/irst index such that $\ell_i(x) < \ell_i(w)$.
Clearly $i$ exists and $\ell_j(x) = \ell_j(w)$ for $j < i$. Let~$i_0'$ be the maximal row index such that $\ell_{i_0'}(w) =
\ell_i(w)$. In other words, \[ i_0' = \max \{i \leq  h\leq n \mid
\ell_{h}(w) = \ell_i(w)\}.\]
Let $j_0$ be the minimal row index such that $\ell_{j_0}(w) \leq \ell_{i_0}(w) - 2$. Notice that $j_0$
exists if $x$ exists: as $\ell_i(x) < \ell_i(w)$, there must be a
row with index $h > i$ such that $\ell_h(x) > \ell_h(w)$. As
$\ell_i(x) \geq \ell_h(x)$, it follows $\ell_i(w) > \ell_h(w) + 1$.
Finally, let
\[ i_0 = \max \{ i' \geq i_0' \mid \ell_{i'}(w)\geq \ell_{j_0}(w) + 2\}.\]
Clearly, $Q = (i_0,j_0)$ is among the
imaginary patterns constructed above, i.e.\ $(i_0,j_0) = (i_r,j_r)$
for some $r$. By construction it is also clear, that $w_Q$ is
$\calP$-stable. Using (\ref{eq:CRITERION}), it follows that $x \leq
w_Q$, and therefore $x = w_Q$.
\end{proof}

\begin{figure}[t]
\centerline{\includegraphics{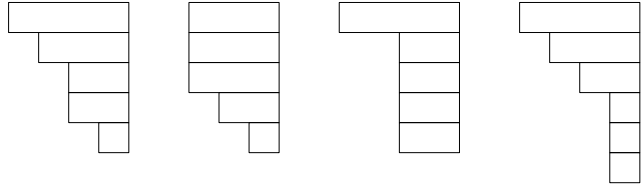}}
\caption{An example of a $\calP$-stable $w$ (left), and the three maximal elements of the singular locus of~$X(w)$. Here $n = 6$, $s = 2$, and the corresponding imaginary patterns are $(1,3)$, $(2,5)$, and $(4,6)$.}\label{fig:PSTABLE} \end{figure}

\begin{rem}As mentioned, any $\calP$-stable $w$ has singularities
unless $w = e$: if $w \neq e$, then $\ell_1(w) > \ell_n(w) + 1$, and
therefore $(1,n)$ is always an imaginary pattern for $w$. In
addition it follows from the description of the $w_{P_r}$ that
\emph{any} $\calP$-stable $w' < w$ is singular. Consequently, the
regular locus of $X(w)$ is just the open orbit $\calP w$. This is well known; see \cite[Theorem~0.1]{em} for a~more general result. It is also shown in~\cite{mov}, where in addition the explicit types of the maximal singularities are described.

In the case $s = 1$, (hence $d = n$),  $X(\ws)$ contains the nilpotent cone, i.e.\ the cone of
nilpotent matrices, as an open af\/f\/ine $\hT$-stable neighborhood $U$
of $e$ (cf.\ Lusztig's isomorphism). In this setting, for any
$\calP$-stable $X(w)$, the intersection $X(w) \cap U$ is the closure of a
nilpotent orbit in~$U$. For these,  again, it is well known that
they are singular along all smaller orbits (\cite{kp}).
\end{rem}

\subsection[$w$ consisting of one string]{$\boldsymbol{w}$ consisting of one string}

The second class of Schubert varieties we will be considering now are those, where the ``relevant'' part of $w$  consists of one string (relevant meaning the part where $w$ dif\/fers from $e$).
In what follows let $e = q_1 < q_2 < \cdots $.

\begin{defn} Let $w \in I^u$. Then $w$ \emph{consist of one string
with critical index $c$} if  for all $k > c$ we
have $w_k = q_k = d(n-1) + k$ (and $c$ is minimal with this
property), and $w_1,w_2,\dots,w_c$ are all congruent mod $n$, i.e.\
there is $j$ such that $w_i \in S_j(w)$ for $1 \leq i \leq c$.
\end{defn}
Examples of such $w$ are of course $e$ (with critical index $c = 0$)
and $\ws$ (resp. $\wsn$) (with critical index $c = d$). More generally for $c =
1,\dots,d$ let $\kappa^c$ be def\/ined as follows
\[ \kappa^c_i = \begin{cases}
(d-c)(n-1) + (i-1)n + 1, & i \leq c,
\\ d(n-1) + i, & i > c.\end{cases}\]
Then $\kappa^d = \ws$. By def\/inition $\kappa^c$ consists of
one string with critical index~$c$. And indeed it is the maximal
such element of $I^u$. In fact, it is the maximum of all elements $w$ in $I^u$ for which $w_i = e_i$ whenever $i > c$ (that is, $w$ need not consist of one string; see the next lemma).

\begin{figure}[t]
\centerline{\includegraphics{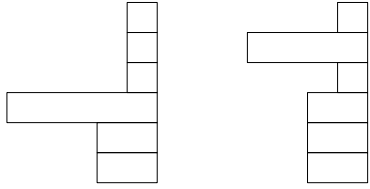}}
\caption{Two elements consisting of one string with critical index $c = 3$ ($n = 6$, $s = 2$). $\kappa^3$ is shown on the left.}\label{fig:ONE_STRING}
\end{figure}

The main reason, why the $w$ consisting of
one string are handled easily, is the following:
\begin{lem}\label{lem:ONE_STRING_BRUHAT}
Let $w \in I^u$ consist of one string with critical index $c$. For
any $x \in I^u$, we have $x \leq w$ if and only if $x_c \geq w_c$,
and $x_i = w_i$ for $i > c$.
\end{lem}
\begin{proof}
We may assume that $c > 1$. The only if part being clear, suppose
$x_i = w_i$ for all $i > c$ and $x_c \geq w_c$. Notice that for $i = 0,1,\dots,c-1$,
\[ w_{c-i} = w_c - in \] since $w$ consists of one string. On the other hand, for any $v \in I^u$, we always have $v_{i+1} - v_i \leq n$, and hence \[ x_{c-i} \geq x_c - in \geq w_c - in = w_{c-i}, \]
and $x \leq w$ follows.
\end{proof}

We will need the following technical criterion below.
\begin{lem}\label{lem:CRIT_ONE_STRING}
Let $w \in I^u$ and suppose $w$ coincides with $e$ at positions $> r$. Then $w$ consists of one string with critical index $c$ if and only if $w_{c-1} \equiv w_c \mod n$ \tu{(}assuming that $c>1$\tu{)}.
\end{lem}

\begin{proof}
The only if part is clear. So suppose $w_{c-1} \equiv w_c \mod n$. Then of course $w_c = w_{c-1} + n$. For any $0 < i < c-1$,  there is $h$ such that $w_{c-1} < w_i + hn \leq w_c = w_{c-1} + n$. But $w_{i} + hn$ belongs to $w$. Hence $w{i} + hn = w_c$.
\end{proof}

Since the Bruhat--Chevalley order is that simple for the case of $w$
consisting of one string, the restriction of $s_\halpha$ being small
in Lemma~\ref{lem:SMALL_REFLECTION} is unnecessary.

\begin{lem}\label{lem:ONE_STRING_REFLECTION}Let $w$ consist of one string. Keeping the notation of
Lemma~{\rm \ref{lem:SMALL_REFLECTION}}, for any $\halpha$ such that
$\xi_\halpha$ is defined at $x \leq w$, we have:
\[ s_\halpha x \leq w \Leftrightarrow r_{i + h + kn,j + kn} x \leq w \text{ for all } k \text{ such  that } 0 \leq k \leq k_0.\]
\end{lem}
\begin{proof} Let us denote $r_{i + h + kn, j+ kn}$ by $r_k$.
Clearly $r_k x \leq w$ whenever $s_\halpha x \leq w$ as
this direction of the assertion holds for any $w$ ($r_k x \leq
s_\halpha x$ if $s_\halpha x > x$); thus, we may assume $r_k x \leq
w$ for all $k$, and it is also safe to assume that $s_\halpha x >
x$. Let $c>0$ be the critical index of $w$. By Lemma~\ref{lem:ONE_STRING_BRUHAT}
we have to show that $(s_\halpha x)_l \geq w_l$ for $l \geq c$.

As $r_k x \leq w$, it is clear that $(r_kx)_l = x_l =
w_l$ for $l> c$, which easily implies that $r_k(x_l) = x_l$ and $(s_\halpha x)_l = w_l$ for
$l > r$. Thus, the only problem might arise if $(s_\halpha x)_c \neq x_c$.
In this case, there is some $k$ such that $(r_k x)_c < x_c$. Thus, $r_k$ moves $x_c$. It follows that $k = r_0$.
Then, if $(r_k x)_c = r_k(x_c)$ we are done, for in this case $(s_\halpha x)_c = r_k(x_c)
\geq w_c$ by assumption. The remaining case is $(r_k x)_c = x_{c-1}$
(equivalent to $r_k (x_c) < x_{c-1}$). Now, if $x_{c-1}$ and $x_c$
are not congruent mod $n$, we are again done, for then $s_\halpha
(x_{c-1}) = x_{c-1} = (s_\halpha x)_c \geq w_c$ because $r_kx \leq w$.

Finally, if
$x_{c-1}$ and $x_c$ are congruent mod $n$, then $x$ itself consists of one string by Lemma~\ref{lem:CRIT_ONE_STRING}.
But if $x$ itself consists of one string, then $(s_\halpha x)_c =
s_\halpha (x_c) = r_k(x_c)$. If $r_k(x_c) < w_c$, then $(r_kx)_c =
x_{c-1} = x_c - n < w_c$ -- a contradiction ($w_c + n$ occurs in both,
$\vert x\vert$ and $\vert w \vert$ at positions strictly bigger than
$c$). Hence, $(s_\halpha x)_c = r_k(x_c) \geq w_c$.
\end{proof}
\begin{cor}\label{cor:ONE_STRING_REAL}
Suppose $\xi_\halpha$ is defined at $x \leq w \in I^u$ where $w$
consists of one string. Then $\xi_\halpha \in T_x(X(w))$ if and only
if $s_\halpha x \leq w$. In particular, \[ T_x(X(w))_\real =
TE(X(w),x).\]
\end{cor}
\begin{proof}
Suppose $\halpha < 0$. Then $s_\halpha x > x$ by Lemma~\ref{lem:REFLECTION_TANGENT}.  If $\xi_\halpha \in
T_x(X(w))$, then $\xi_\halpha \in T_x(Y(w))$, and therefore (using
the notation of Lemma~\ref{lem:SMALL_REFLECTION} and the proof of Lemma~\ref{lem:ONE_STRING_REFLECTION}), $r_k x \leq w$ for all $0 \leq k
\leq k_0$. Lemma~\ref{lem:ONE_STRING_REFLECTION} now gives $s_\halpha x
\leq w$. All other cases are immediate.
\end{proof}
In the case of imaginary tangents, a similar result holds:
\begin{lem}\label{lem:ONE_STRING_IMAGINARY}
Let $w \in I^u$ consist of one string with critical index $c$.  If
for any $x \leq w$, $S(h,x,w)$ \tu{(}cf.\ \eqref{eq:SHXW}\tu{)} contains $i \neq j$, say, then
$\xi_{(ij),h} \in T_x(X(w))$. In particular, $x$ is singular.
\end{lem}
\begin{proof}
Without loss of generality, $h_i(x) < h_j(x)$. Let $m$ be the unique
nonnegative integer such that $h_j(x) - h_i(x) = (j-i) + mn$. Let
$\hbeta = (ij) - m\delta$. Then $\hbeta < 0$ by construction.
Moreover, $s_\hbeta x = x$ because $s_\hbeta(h_j(x)) = h_i(x)$. Furthermore, $s_{\hbeta - h\delta}x \leq w$. To see this, notice
that as $S(x,h,w)$ contains two or more elements, $x$ cannot consist
of one string with critical index greater or equal $c$. It is clear
that $s_{\hbeta-h\delta} x \in I^u$ (because
$\ell_j(x) \geq h$ and $h_i(x) - hn > 0$), and is obtained from $x$ by increasing $\ell_i(x)$ and decreasing $\ell_j(x)$ by $h$. By hypothesis, $w$ consists of one
string, so we have to see that $(s_{\hbeta-h\delta}x)_k = w_k$ for
$k > c$, and $(s_{\hbeta - h \delta }x)_c \geq w_c$. The f\/irst
assertion is clear, because $s_{\hbeta-h\delta}$ changes only those
elements of $x$  which are also changed by one of the
$r_{h_i(x)+(k-h)n,h_i(x) + kn}$ or $r_{h_j(x) + (k-h)n,h_j(x) + kn}$, and these do not change the entries $x_l$ for $l > c$.

There are two possibilities: Either
$s_{\hbeta -h\delta}$ does not change $x_c$ and we are done, or
$(s_{\hbeta-h\delta}x)_c$ is $x_{c-1}$ with $x_{c-1} \not \in
S_j(x)$: $s_{\hbeta -h\delta}(x_c) \neq x_c$, therefore $x_c \in
S_j(x)$ is the largest element changed, i.e.\ $x_c = h_j(x) + (h-1)
n$. As $x$ does not consist of one string with critical index $c$,
$x_{c-1}\not \in S_j(x)$ (Lemma~\ref{lem:CRIT_ONE_STRING}), and $x_{c-1} = (r_{h_j(x)-n,h_j(x) +
(h-1)n}x)_c$ (recall that $x_{c-1} > x_c - n$). As $j \in S(h,x,w)$, this last statement means $x_{c-1} \geq w_c$ and so $s_{\hbeta - h\delta}x \leq w$.

But now we conclude, as in the proof of Lemma~\ref{lem:IMAGINARY_SING} that $\xi_{(ij),h} = \pm [\xi_{-\hbeta},\xi_{\hbeta - h\delta}]$ is tangent to $X(w)$ at $x$.
\end{proof}
Summarizing, we have obtained:
\begin{thm}\label{thm:ONE_STRING_TANGENT}Suppose $w \in I^u$ consists of one string. Then for all
$x \leq w$ we have $T_x(X(w)) = T_x(X(\ws)) \cap T_x(Y(w))$. In
particular, in $T_x(Y(w))^u$, $T_x(X(w))$ is given by the trace
relations.
\end{thm}
\begin{proof}
Lemma~\ref{lem:ONE_STRING_REFLECTION} in particular says that if $\xi_{\halpha}$ is def\/ined at $x$ and is contained in $T_x(Y(w))$, then it is contained in $T_x(X(w))$. By Lemma~\ref{lem:REFLECTION_TANGENT}, this means that
\[ T_x(X(w))_{\real} = T_x(Y(w))^u_{\real} = T_x(Y(w)) \cap T_x(X(\ws))_{\real}. \]

Regarding imaginary roots, Lemma~\ref{lem:ONE_STRING_IMAGINARY} immediately implies that the imaginary part of $T_x(X(\ws))$ is spanned by all $\xi_{(ij),h}$ for which $i \neq j \in S(h,x)$. In other words, $T_x(X(\ws))_{-h\delta}$, as a subspace of $T_x(G(d,V))_{-h\delta}$, is cut out by the trace relation for $h$ (cf.\ Lemma~\ref{lem:TRACE_RELATION}).
Applying Lemma~\ref{lem:ONE_STRING_IMAGINARY} in the case of arbitrary $w \in I^u$, this, combined with Lemma~\ref{lem:IMAGINARY_TANGENTS}, means that for each $h$, $\xi_{(ij),h} \in T_x(X(w))$ if and only if $\xi_{(ij),h} \in T_x(X(\ws))_{-h\delta} \cap T_x(Y(w))$.
This completes the proof.
\end{proof}

\begin{rem}\label{rem:ONESTRING_SINGLOCUS}
Theorem~\ref{thm:ONE_STRING_TANGENT} together with Corollary~\ref{cor:ONE_STRING_REAL} immediately imply that for $w$ consisting of one string, $X(w)$ is smooth at $x \leq w$ if and only if the number of $\hT$-stable curves in $E(X(w),x)$ equals the dimension of $X(w)$, and if there is no imaginary tangent in $T_x(X(w))$. This is consistent (though stronger) with the observation in Remark~\ref{rem:Peterson}, which on the other hand applies to arbitrary $w$ not just those consisting of one string.
\end{rem}

Let $w$ consist of one string with critical index $c > 2$. If $w =
\kappa^c$ let $P = ([w_1],[w_1+1])$, otherwise put $P =
([w_1],[w_{c+1}-n])$, where we put $w_{c+1} = d(n-1) + c + 1$ in
case $c = d$. Clearly~$P$ is an imaginary pattern for $w$, and we
def\/ine $\varphi(w) = w_P$.

Thus, if $w \neq \kappa^r$, $\varphi(w)$ is
obtained from $w$ by replacing $w_1$ with $w_{c+1}-n$; and
$\varphi(\kappa^c)$ is obtained by replacing $\kappa^c_1$ with
$\kappa^r_{c} + 1$.

\begin{figure}[t]
\centerline{\includegraphics{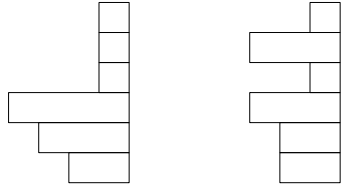}}
\caption{$\varphi(w)$ for the elements consisting of one string shown in Fig.~\ref{fig:ONE_STRING}. $\varphi(\kappa^3)$ is left.}\label{fig:PHI_W}
\end{figure}

It is clear that $\varphi(w)$ is singular. But
in fact we have
\begin{thm}\label{thm:ONE_STRING_SING}Let $w$ consist of one string with critical index $c >
2$. Then the singular locus of $X(w)$ is $X(\varphi(w))$.

In particular, the singular locus of $X(w)$ has exactly codimension two.
\end{thm}
For the proof we will need:
\begin{thm}\label{thm:RATIONALLY_SMOOTH}
Let $x \leq \kappa^c$ for some $c > 0$. Then the set
\[ E(\kappa^c,x) = \{ s \in \hW \mid x \neq sx \leq \kappa^c\} \]
has precisely $c(n-1) = \dim X(\kappa^c)$ elements. In particular
$\dim T_x(X(\kappa^c))_\real = \dim TE(X(\kappa^c),x)$ $= \dim
X(\kappa^c)$.
\end{thm}
\begin{proof}
The assertions on the dimension of $T_x(X(\kappa^c))$ and
$TE(X(\kappa^c),x)$ are immediate consequences of the f\/irst
assertion and Corollary~\ref{cor:ONE_STRING_REAL}.

We proceed by descending induction on $c$. Suppose $c = d$. Then
$\kappa^c = \ws$, and we simply have to count all ref\/lections
$s_\halpha$, def\/ined at $x$, with $s_\halpha x \neq x$. It will be
convenient to use the notation introduced in the beginning: we write
$x  = \tau^{(c_1,c_2,\dots,c_n)}$, then $h_i(x) = i + c_in$. We have
to count the def\/ined $\xi_\halpha$s. $\xi_\halpha$ may map a given
$e_{h_i(x)}$ to any $e_j$, provided $j \not \equiv i \mod n$, $j
\not \in \vert x \vert$, and $h_{[j]}(x) - j \leq n\ell_i(x)$ (the
latter because $\xi_\halpha$ is $\tau$-equivariant).  By
construction it is clear that also $h_{[j]}(x) = [j] + c_{[j]}n$ thus
there are at most $c_{[j]}$ such maps sending $e_{h_i(x)}$ to an
element congruent to $[j]$.

If $E = \vert E(\kappa^c,x)\vert$ we  conclude
\[ E = \sum_{i = 1}^n \sum_{\substack{1\leq j \leq n\\j \neq i}}
\min\{c_j,\ell_i(x)\}.\] Of course, $d = \sum_i \ell_i(x)$, and $c_j
= d - \ell_j(x) = \sum_{i\neq j} \ell_i(x)$. Thus, $c_j \geq
\ell_i(x)$ for all $i \neq j$, which implies that
\[ E = \sum_i \sum_{j\neq i} \ell_i(x) = d(n-1).\]This settles the case $c = d$: we already know  $TE(X(\ws),x)$ contains all $\xi_\halpha$ that are def\/ined (cf.\ Lemma~\ref{lem:REFLECTION_TANGENT}). We f\/ind that $\dim TE(X(\ws),x) = d(n-1)$ for all $x \leq \ws$. This applies in particular to $x = \ws$ showing that $\dim X(\ws) = d(n-1)$.

Now suppose $c < d$ and the assertion
holds for $c' > c$. The ref\/lections which we have to count are
precisely those which f\/ix
$\kappa^c_{c+1},\kappa^c_{c+2},\dots,\kappa_d^c$: $s_\halpha x \leq
\kappa^c$ if and only if
\begin{gather}
 \label{eq:FIXED_POSITIONS}
(s_\halpha x)_{k} = \kappa^c_k
\end{gather} for
$k>c$ and $(s_\halpha x)_c \geq \kappa^c_{c}$. This last condition
is void because $\kappa^c_{c} = \kappa_{c+1}^c - n$.

By induction we may assume that the number of ref\/lections for which \eqref{eq:FIXED_POSITIONS} holds for $k > c+1$
equals $\dim
X(\kappa^{c+1}) = (c+1)(n-1)$ (note that $\kappa^c_k = \kappa^{c+1}_k$ for these $k$). Let $E'$ be the set of these ref\/lections. For each $i \in \{1,2,\dots,n\}$ which is not congruent
$\kappa^c_{c+1}$ mod $n$, there is exactly one ref\/lection $s$ in $E'$
for which $(sx)_{c+1} \neq x_{c+1} = \kappa_{c+1}^c$: namely, the one moving the
entire $n$-string in $x$ through $x_{c+1}$ to the string through
$i$. To be precise, if the congruence class of $x_{c+1}$ is $j$,
then $s$ replaces $x$ by removing all elements in
$S_j(x)$ up until (including) $x_{c+1}$ and replacing them by adding the same
number of elements to the string through $i$. This is always
possible (as $\ell_j(x) \leq c_i$). Of course, $s$ is an up-exchange, because $x_{c+1}$ is moved to a smaller number. Also, $(sx)_{c+1} < x_{c+1}$, and therefore $sx \nleq \kappa^c$.

Thus $E = c(n-1)$. As $x$ was
chosen arbitrary, this in particular says that
$E(\kappa^c,\kappa^c)$ has $c(n-1)$ elements. $\kappa^c$ is a smooth
point of $X(\kappa^c)$, so here we know that the number of curves is equal to the  dimension
of $X(\kappa^c)$, and therefore $\dim X(\kappa^c) = c(n-1)$.
\end{proof}
\begin{rem}An immediate consequence of this theorem is the fact that
$X(\kappa^c)$ is rationally smooth (cf.\ \cite{CP}), which of
course is well known, at least in the case of $c = d$, and also shown in~\cite{BM}.

Furthermore, as is easily seen, $T_e(X(w))_\real =
T_e(X(\kappa^c))_\real$ whenever $w$ consists of one string with
critical index $c \geq 2$. Therefore, for such a $w$ dif\/ferent from $\kappa^c$, $X(w)$
cannot be globally rationally smooth, as $\dim X(w) < \dim X(\kappa^c)$.
\end{rem}

Returning to the situation of an arbitrary $w$ consisting of one
string, we have:
\begin{lem}\label{lem:ONE_STRING_ONE_STRING}Suppose $x \leq w$ both consist of one string with
critical index $c > 0$. Then $x$ is a~regular point of $X(w)$.
\end{lem}
\begin{proof}
Let $E = \{ s_\halpha \mid x \neq s_\halpha x \leq \kappa^c\}$. By
Theorem~\ref{thm:RATIONALLY_SMOOTH}, we know that $\vert E \vert =
c(n-1)$. Furthermore the codimension $c$ of $X(w)$ in $X(\kappa^r)$
is precisely $w_1 - \kappa_1^c = w_c - \kappa^c_c$ (by Lemma~\ref{lem:CODIM}). Now
$x$ consists of one string with critical index $c$, and therefore
any $s_\halpha \in E$ satisf\/ies $s_\halpha x \leq w$ if and only if
$s_\halpha(x_c) = (s_\halpha x)_c \geq w_c$. There are precisely $c$
elements $s_\halpha$ of $E$ such that $\kappa^c_c \geq
s_\halpha(x_c) > w_c$: if $s_\halpha x > x$ and $s_\halpha\in E$,
then $s_\halpha x$ consists of one string, so $s_\halpha x$ is
uniquely determined by~$s_\halpha(x_c)$. As a consequence $\vert
E(X(w),x) \vert = \dim X(\kappa^c) - c = \dim X(w)$.

It remains to show that $T_x(X(w))$ has no imaginary weight.
However, as $T_x(X(w)) \subset T_x(X(\kappa^c))$ this may be checked
in case $w = \kappa^c$. So suppose $\xi_{i,h}$ is def\/ined at $x$ for
some~$i$,~$h$, and contained in $T_x(Y(\kappa^c))$. Then
$r_{h_i(x)-hn,h_i(x)}x \leq \kappa^c$. In particular, $h_i(x)$ must
appear before the critical index $c$ of $x$ and thus $h_i(x) = x_1$;
consequently, no other $\xi_{j,h}$ is contained in
$T_x(Y(\kappa^c))$ for this given $h$. The trace condition now kills
$\xi_{i,h}$ (in fact, as $h_i(x) = x_1$, and $x_1 - \kappa^c_1 < n$,
even~$\xi_{i,h}$ is not tangent to $T_x(Y(\kappa^c))$.
\end{proof}
\begin{rem}\label{rem:ONE_STRING_INDEX_1} Notice that this already shows the (elementary) fact that
$X(w)$ is globally nonsingular for all $w$ consisting of one string with
critical index $c = 1$.
\end{rem}

\begin{proof}[Proof of Theorem~\ref{thm:ONE_STRING_SING}]
We have to show that any $x \leq w$ which is not below $\varphi(w)$
is a~regular point of $X(w)$. By
Lemma~\ref{lem:ONE_STRING_ONE_STRING} it suf\/f\/ices to treat the
case when $x$ does not consist of one string with critical index $c$
itself.

Let $l = \vert \{ x_i \mid w_{c+1}-n \leq x_i < w_{c+1}\}\vert$ be
the number of entries of $x$ between $w_{c+1} - n$ and $w_{c+1}$. If
$l = 1$, then $x = \kappa^{c-1}$ since we excluded the remaining
possibility that $x$ consists of one string with critical index $c$.

Suppose f\/irst that $x_1 \geq w_2$. Then $x_{2} \geq
w_3,\dots,x_{c-2} \geq w_{c-1}$. Otherwise $x_{c-2} < w_{c-1}$ (as
all $w_i$ lie on one string), and thus $x_h = x_{c-1}$ is the only
entry of $x$ with $w_{c-1} \leq x_h < w_c$. Consequently the entries
$x_1,x_2,\dots,x_{c-1}$ are all congruent mod $n$, and therefore
$x_1 < w_{c-1}-(c-3)n = w_2$, a contradiction. But now $x \leq
\varphi(w)$, since $x_{c-1} \geq w_{c+1}-n$: If $x_{c-1} <
w_{c+1}-n$, then $x_c = x_{c-1} + n < w_{c+1}$, and thus $l = 1$,
and therefore $x = \kappa^{c-1}$. Again, this contradicts $x_1 \geq
w_2$. Summarizing, if $x_1 \geq w_2$ then $x \leq \varphi(w)$
(obviously this is only-if as well).

It remains to treat the case when $x_1 < w_2$. In this case $l = 2$,
or $x = \kappa^{c-1}$. In the second case, the only ref\/lections $s$
with $\kappa^{c-1}< s\kappa^{c-1} \leq w$ are precisely the
up-exchanges of $i = [\kappa^{c-1}_1]$ with a number between $i$ and
$[w_1]$, giving a total of $\kappa^{c-1}_1 - w_1$ which is the
codimension of $x = \kappa^{c-1}$ in $X(w)$. In the other case ($l =
2$), we have $x_1 \equiv x_2 \equiv \dots \equiv x_{c-1} \mod n$ but
$x_1 \not \equiv x_c \mod n$. Thus, the only possibilities for a
ref\/lection are to move $x_{c-1}$ to any integer between $w_{c+1}-n$
and $x_{c-1}$, or $x_c$ to an integer between $w_c$ and $x_c$, or
f\/inally, $x_{c-1}$ to $x_{c}-n$ (one might think one could move
$x_{c-1}$ also to values between $w_{c-1}$  and $w_{c+1}-n$, but
that is not possible in general since such a ref\/lection would have
to move $x_{c-1} + n \geq w_{c+1}$ as well. Counting these
ref\/lections gives $(x_{c-1} - w_{c+1}-n) + (x_c - w_c) + 1$. But
this is just the codimension of $X(x)$ in $X(w)$ ($x$ is obtained from $w$ by f\/irst
down-exchanging $[w_1]$ and $[x_1]$, and then down-exchanging
$[w_c]$ and $[x_c]$). Concluding it follows that $E(X(w),x)$ has the
minimal number of elements possible, namely $\dim X(w)$. On the
other hand, there is no imaginary tangent at $x$: if $l = 1$ the
only possibility for such a tangent is $\xi_{[x_1],1}$ but $x_1-n <
w_1$. If $l = 2$ there is at most one other possibility
(corresponding to $[x_c]$) but this is killed by the trace relation
(or the remark that $x_c-n < w_c$ and $x_{c-1} < w_c$).
Theorem~\ref{thm:ONE_STRING_TANGENT} now gives the result.

Finally, that $X(\varphi(w))$ has codimension two in $X(w)$ follows from the fact that here $\varphi(w) = stw$ where $s,t$ are small ref\/lections and $stw < tw< w$ are codimension one steps (cf.\ Lem\-ma~\ref{lem:CODIM}).~\end{proof}

\begin{rem}
In view of Remark~\ref{rem:ONE_STRING_INDEX_1} the only remaining
case is $w$ consisting of one string with critical index $c = 2$.
The only dif\/ference to the case $c > 2$ is that $\varphi(w)$ has to
be def\/ined slightly dif\/ferent: the pattern for such a $w$ is $P =
([w_1],[w_1] + 1)$ (notice that $[w_1] = [w_2]$). All the proofs
above go through when this is kept in mind appropriately. The main
dif\/ference now is that $X(\varphi(w))$ does not have codimension two
as in the case of critical index $c > 2$.
\end{rem}

\begin{rem}\label{rem:SMOOTH}\looseness=-1
Finally, let us conclude with a short remark on the smooth Schubert varieties:
In~\cite{BM} Billey and Mitchell completely classify all smooth Schubert varieties in af\/f\/ine Grassmannians for all types. They all are (closed) orbits for certain parabolic ind-subgroups of $\mathcal G$. In type $A$, their result means (in our notation) that the smooth $X(w)$ are precisely those  where $L(w) = (l_1,l_2,\dots,l_{2n})$ has the following form (we list only $l_1,l_2,\dots,l_n$):
There is a pair of integers $1 \leq p < q \leq n - p$ such that $l_1 = l_2 = \dots = l_p = 1$ and $l_q = l_{q +1} = \dots = l_{q + p} = 2$ and $l_i = 1$ for all remaining $i$ between $p+1$ and $n$.
Alternatively (somewhat dual to this construction) there is a second ``family'' given by integers $1 \leq p < q \leq n - p$ with $l_{n-1} = l_{n-2} = \dots  = l_{n-p} = 2$ and $l_{n-q} = l_{n-q-1} = \dots = l_{n-q - p} = 0$.
In particular, all of them are below $w^1$ (i.e.\ $s  = 1$ and $d = n$).

This is consistent with our discussion: if $X(w)$ is smooth it cannot admit any imaginary pattern and hence $l_i - l_j$ cannot be strictly greater than $1$ if $i < j$ and $2$ if $i > j$. This already shows that $w \leq w^1$, and hence $l_i \in \{ 0,1,2\}$ for $i \leq n$. Not allowing any real pattern then asserts that we are in one of the two cases listed. Indeed, it follows that all the $i$ with $l_i = 2$ are strictly larger than all the $i$ with $l_i = 0$. Also. there cannot be any ``gap'' between $i$ and $j > i$ for which $l_i = l_j = 2$; similarly, all the $l_i$ with $l_i = 0$ must be consecutive as well.
Assuming $w \neq e$, if $l_1 \neq 0$, then $l_n = 2$ for otherwise there is a real pattern $(i,n,n+1,j)$ where $[j] < i$ and $l_j = 1$.

It seems plausible that one could show the smoothness of these Schubert varieties also from our discussion (maybe using Remark~\ref{rem:Peterson}) by arguing that the classical Schubert varieties $Y(w)$ won't contain any imaginary tangents for all these varieties.

In \cite{BM} the rationally smooth Schubert varieties are also classif\/ied. There are two types (other than the smooth ones): what the authors call ``spiral'' corresponds to our $\kappa^c$ (and a dual version, related by an automorphism) and we showed above that this is indeed rationally smooth.
The second type is ``chains'', that is, Schubert varieties, that contain exactly one f\/lag of subvarieties $X(e) \subsetneq X(w_1) \subsetneq X(w_2) \subsetneq \cdots \subsetneq X(w_n) = X(w)$ where $n = \dim X(w)$. According to \cite{BM}, except for the case $n = 1$, these are all smooth.

Hence the only rationally smooth Schubert varieties are the smooth ones and (essentially) those of the form $X(\kappa^c)$. As one referee suggested, it is tempting to conjecture that this means all $w$ except $w = \kappa^c$ or those for which $X(w)$ is smooth, should admit a real pattern of some sort.
While this may be true, consider the following example: $n = d = 4$, so $s = 1$, and $w = s_{(12)}\ws$. Then $w$ consists of one string with critical index $d$. It is not rationally smooth since it is strictly smaller than $\kappa^d = \ws$.
However the only real pattern supported by $w$ is degenerate: it is an exceptional real pattern of the f\/irst kind, namely $P = (2,5,6,7)$ ($L(w) = (0,4,0,0,1,5,1,1)$). The corresponding singular point $w_P$ is def\/ined as $L(w_P) = (0,3,1,0,1,4,2,1)$. Note that $w_P < \varphi(w)$ and indeed, $w_P$ admits four up-exchanges whereas the codimension of $X(w_P)$ in~$X(w)$ is only three. $\varphi(w)$ on the other hand is still a rationally smooth point of $X(w)$.
\end{rem}

\subsection*{Acknowledgments}

The f\/irst author was supported by the Swiss National Science Foundation, and
partially by an NSERC Discovery Grant. The second author was
supported by NSF grant DMS-0652386 and Northeastern University
RSDF 07--08.
The f\/irst author would like to thank the Swiss National Science Foundation for making possible his stay at Northeastern University, during which most of this work has been done.
We would like to thank the referees for their careful reading and their many valuable suggestions.

\pdfbookmark[1]{References}{ref}
\LastPageEnding

\end{document}